\documentclass[12pt]{amsart}

\usepackage{amsmath}
\usepackage{amssymb}
\usepackage{amsthm}
\usepackage{xcolor}
\usepackage{bm}
\usepackage{graphicx}
\usepackage{pifont}
\usepackage{mathrsfs}
\usepackage{epsfig,verbatim}
\usepackage[pagewise]{lineno}\linenumbers

\newcommand\be{\begin{equation}}
\newcommand\ee{\end{equation}}
\newcommand\bea{\begin{eqnarray}}
\newcommand\eea{\end{eqnarray}}
\newcommand\beaa{\begin{eqnarray*}}
\newcommand\eeaa{\end{eqnarray*}}
\newcommand\beba{\begin{equation}\left\{\begin{array}{rcl}}
\newcommand\eeba{\end{array}\right.\end{equation}}
\newcommand\bebaa{\begin{equation*}\left\{\begin{array}{rcl}}
\newcommand\eebaa{\end{array}\right.\end{equation*}}

\newcommand\bR{{\mathbb{R}}}

\newcommand\bN{{\mathbb{N}}}
\newcommand\bZ{{\mathbb{Z}}}
\newcommand\cA{{\mathcal{A}}}

\newcommand\cH{{\mathcal{H}}}
\newcommand\cL{{\mathcal{L}}}
\newcommand\cR{{\mathcal{R}}}

\newcommand\cN{{\mathcal{N}}}

\newcommand\cX{{\mathcal{X}}}
\newcommand\cY{{\mathcal{Y}}}

\newcommand{\dsum}{\displaystyle \sum}

\newcommand{\dlim}{\displaystyle \lim}
\newcommand{\dmax}{\displaystyle \max}

\newtheorem{theorem}{Theorem}[section]
\newtheorem{lemma}[theorem]{Lemma}
\newtheorem{corollary}[theorem]{Corollary}

\newtheorem{definition}[theorem]{Definition}

\newtheorem{proposition}[theorem]{Proposition}
\newtheorem{assumption}[theorem]{Assumption}

\numberwithin{equation}{section}
\def\theequation{\arabic{section}.\arabic{equation}}

\begin{document}
\nolinenumbers

\title[Turing bifurcation]{Turing bifurcation, criticality, and multistability in a nonlocal advection-reaction-diffusion equation}

\author[Y. Nagayama]{Yui Nagayama}
\address{Department of Mathematics, Hokkaido University, Hokkaido, 060-0810, Japan}
\email{nagayama.yui.z5@elms.hokudai.ac.jp}

\author[M. Nagayama]{Masaharu Nagayama}
\address{Research Center of Mathematics for Social Creativity, Research Institute for Electronic Science, Hokkaido University, Hokkaido, 060-0812, Japan}
\email{nagayama@es.hokudai.ac.jp}

\author[H. Ishii]{Hiroshi Ishii}
\address{Research Center of Mathematics for Social Creativity, Research Institute for Electronic Science, Hokkaido University, Hokkaido, 060-0812, Japan}
\email{hiroshi.ishii@es.hokudai.ac.jp}

\thanks{Date: \today. Corresponding author: Hiroshi Ishii}



\begin{abstract}
This paper considers a nonlocal advection--reaction--diffusion equation with logistic growth. 
We investigate how the interaction kernel affects the criticality of Turing bifurcations and the associated multistability. 
We first prove that nonhomogeneous stationary solutions bifurcate from the positive homogeneous stationary solution.
Center manifold reduction yields an amplitude equation near the bifurcation point. 
This equation provides criteria for supercritical and subcritical bifurcations in terms of the Fourier transform of the interaction kernel. 
We also derive a sufficient condition for subcriticality. 
This condition can be satisfied by sufficiently localized attractive interactions. Based on these analytical results, we investigate multistability through numerical bifurcation analysis. 
The numerical results show that the stable homogeneous stationary solution can coexist with stable nonhomogeneous stationary solutions or stable time-periodic solutions.
\end{abstract}

\maketitle
\setcounter{tocdepth}{3}

\section{Introduction}

The spatial organization and pattern formation of biological populations and cell collectives have long been studied using mathematical models.
Among various modeling frameworks, reaction--diffusion systems provide a fundamental approach for describing such phenomena. 
By combining population growth and interactions with random spatial movement, these systems have played a central role in the mathematical study of pattern formation.
Diffusion-driven instability, originally proposed by Turing \cite{T52}, is one of the fundamental  mechanisms underlying pattern formation.
It has provided a theoretical basis for many subsequent studies (e.g., \cite{KM10, Meinhardt, Murray}).
This instability can also be studied from the viewpoint of bifurcation theory.
As a model parameter is varied, nonhomogeneous stationary solutions may bifurcate from a homogeneous stationary solution at the onset of the instability.
Such a bifurcation is commonly referred to as a \emph{Turing bifurcation} \cite{BMS98, BC14, DFM15}.

Experimental studies have shown that cells and organisms do not always move randomly.
They can also sense their environment and other individuals and change their direction of movement accordingly \cite{MT15, Remote2020, YK14}. 
Nonlocal advection models have attracted attention as a framework for describing such directed movement.
Nonlocal advection represents a mechanism in which individuals at a given location sense the distribution of other individuals in the surrounding region and move in response to this information.

Nonlocal advection can be derived from several modeling perspectives, including cell-cell adhesion \cite{APS06, MT15}, reductions of chemotaxis models \cite{IMT, MT25}, and continuum limits of interacting particle systems \cite{G24}.
A representative nonlocal advection equation is given by
\be\label{eq:nae}
    u_t = - [u (K*u)_{x}]_{x}.
\ee
Here, $u(t,x)$ denotes the population density, $K(x)$ is an even interaction kernel, and $K*u$ denotes the convolution integral
\beaa
    (K*u)(t,x) := \int_{\bR} K(x-y)u(t,y)dy.
\eeaa
The interaction kernel $K$ determines the nonlocal interaction, and different choices of $K$ can produce different spatial patterns \cite{CMSTT19, JKMG23}.
In fact, equation \eqref{eq:nae} can be derived as the mean-field approximation of the particle system
\beaa
    \dfrac{d x_l}{dt} =   -\dsum^M_{m =1} K'(x_m-x_l).
\eeaa
This particle system can be interpreted as the gradient flow of the energy
\beaa
    E[x_1,\ldots,x_M] = -\dfrac{1}{2} \dsum^M_{m=1} \dsum^M_{l=1} K(x_{m} - x_{l}).
\eeaa
Therefore, if $K$ is monotonically decreasing for $x>0$, the particles move toward each other, whereas if $K$ is monotonically increasing for $x>0$, they move away from each other.
Thus, the properties of the interaction kernel $K$ determine the direction and strength of the interactions between individuals and strongly affect the structure of solutions and the resulting spatial patterns.


Mathematical models with nonlocal advection have been widely used to describe the dynamics of cell populations. 
Armstrong, Painter, and Sherratt \cite{APS06} proposed a model in which cell--cell adhesion is represented by a nonlocal advection term. 
They showed that Turing-type patterns can arise even in a scalar equation. 
The model was later extended by Murakawa, Togashi, and their collaborators \cite{CMSTT19, MT15}, and these extensions reproduced experimental observations such as cell sorting.
Nonlocal advection models have also been used to study other biological phenomena, including run-and-chase interactions between pigment cells in zebrafish, territorial pattern formation, and collective behavior \cite{JKMG23, PHP24}. 
These studies show that nonlocal advection provides a flexible framework for describing self-organized spatial patterns and complex structures.

In this study, we consider a nonlocal advection--reaction--diffusion equation and examine how the interaction kernel affects the bifurcation structure of spatial patterns.
The equation is posed on the interval $I:=(0,L)$, with $L>0$, subject to periodic boundary conditions, and takes the form
\be\label{eq:main}
    u_t = du_{xx} - k[u(K*u)_x]_x + f(u),
    \qquad t>0,\quad x\in I.
\ee
Here, $u=u(t,x)\in\bR$ represents the population density, $d>0$ is the diffusion coefficient, $k>0$ measures the strength of nonlocal advection, and $f:\bR\to\bR$ is the reaction term.
We assume that the interaction kernel $K\in L^1(\bR)$ is even.
For $L$-periodic functions $u$, we define the convolution by
\beaa
    (K*u)(t,x)
    =
    \sum_{j\in\mathbb{Z}}
    \int_I K(x-y-jL)u(t,y)\,dy.
\eeaa

Equation \eqref{eq:main} and related models have been studied in several settings. 
When $f = 0$, equation~\eqref{eq:main} has been extensively studied as an aggregation--diffusion equation (e.g., \cite{CCY19,G24,PHP24}). 
The equation admits a free-energy structure.
In this setting, a central problem is to understand the long-time behavior of solutions and the structure of stationary solutions.
Many studies have therefore focused on the corresponding stationary problem; see the review \cite{CCY19}.


On the other hand, when the reaction term $f(u)$ is present, several aspects of the equation have been studied, including stationary problems \cite{Ei86,MO86}, propagation phenomena \cite{HH20}, and Turing bifurcations \cite{DFM15,JKMG23}. Ducrot, Fu, and Magal \cite{DFM15} analyzed the linearized equation around a homogeneous stationary solution.
They then applied the Crandall--Rabinowitz bifurcation theorem to prove the existence of a Turing bifurcation in a nonlocal advection--reaction--diffusion equation.
Jewell et al. \cite{JKMG23} studied how pattern formation depends on the spatial dimension. 
Their linear stability analysis and numerical simulations showed that Turing bifurcations can occur.
Their numerical results also suggested the occurrence of subcritical Turing bifurcations.

A subcritical bifurcation can lead to multistability. 
In particular, a stable homogeneous solution may coexist with stable nonhomogeneous solutions.
Such multistability is known to give rise to hysteresis and to different spatial patterns depending on the initial condition (e.g., \cite{BC14,Strogatz}).
Determining the criticality of a Turing bifurcation helps clarify the stability of the emerging patterns. 
It also provides a basis for studying multistability and transitions between different states.
However, it remains unclear how the properties of the interaction kernel $K$ affect the criticality of Turing bifurcations in nonlocal advection--reaction--diffusion equations.


The aim of this paper is to clarify the relation between the properties of the interaction kernel and the criticality of Turing bifurcations.
We also investigate what types of multistability can arise from a subcritical bifurcation.
More specifically, for equation \eqref{eq:main}, we consider the logistic reaction term
\be\label{logistic}
    f(u) = c u(1-u) \quad (c>0),
\ee
which describes growth and saturation.
First, we use the Crandall--Rabinowitz theorem to prove the existence of nonhomogeneous stationary solutions bifurcating from the homogeneous stationary solution.
Next, we derive an amplitude equation near the bifurcation point by center manifold reduction.
The cubic coefficient determines whether the bifurcation is supercritical or subcritical. Its sign is characterized in terms of the Fourier transform of $K$.
We then use numerical bifurcation analysis to trace the branch of nonhomogeneous stationary solutions over a wide parameter range.
In this process, we also investigate time-periodic solutions that arise through secondary bifurcations of the nonhomogeneous stationary solutions.
This enables us to examine multistability involving both stationary and time-dependent patterns.


This paper is organized as follows.
In Section~\ref{sec:sta}, we present the results on bifurcation of stationary solutions and derive an amplitude equation near the bifurcation point.
We also analyze the criticality of the bifurcation and derive conditions for the occurrence of a subcritical Turing bifurcation.
In Section~\ref{sec:num}, we perform numerical bifurcation analysis under these conditions and examine how multistability arises.
Finally, in Section~\ref{sec:dis}, we summarize our results and discuss open problems.


\bigskip
\section{Local bifurcations of stationary solutions}\label{sec:sta}
In this section, we study local bifurcations of stationary solutions to equation \eqref{eq:main} with the logistic reaction term \eqref{logistic}.
The equation has two homogeneous stationary solutions, $u = 0$ and $u = 1$.
Since $u = 0$ is always unstable, we focus on bifurcations from $u=1$, taking $k>0$ as the bifurcation parameter.

First, we apply the Crandall--Rabinowitz theorem to establish the local bifurcation of nonhomogeneous stationary solutions from the homogeneous stationary solution.
We then derive an amplitude equation for the critical mode by means of a center manifold reduction and analyze its higher-order terms to determine the local structure and criticality of the bifurcation.
Finally, we characterize conditions on the interaction kernel under which a subcritical bifurcation occurs.

\subsection{Existence of nonhomogeneous stationary solutions}

Throughout this section, we consider the reaction term $f(u)$ given by \eqref{logistic}.
We first examine the linear stability of $u=1$.
Setting $u = 1 + w$ in equation \eqref{eq:main} gives
\begin{align}
    w_t = dw_{xx} -cw -cw^2 -k(K*w_{xx}) -k[w(K*w)_x]_x =: \mathcal{H}(k, w; d, c). \label{eq:main_expand}
\end{align}
We define the linear and nonlinear parts, respectively, by
\begin{align*}
    \cL w=\mathcal{L}(k,d,c)w &:= dw_{xx} -cw -k(K*w_{xx}), \\
    \cN (w)=\cN(w;k,c) &:= -cw^2 -k[w(K*w)_x]_x.
\end{align*}
Both operators are defined as mappings from $H^2_{\mathrm{per}}(I)$ into $L_{\mathrm{per}}^2(I)$.
Here, $X_{\mathrm{per}}(I)$ denotes the subspace of $X(I)$ consisting of all $L$-periodic functions.
We note that $\cN$ is well-defined since $H^1_{\mathrm{per}}(I) \subset L^{\infty}_{\mathrm{per}}(I)$ by Morrey's inequality.

We expand $w$ in a Fourier series as
\begin{align}
    w(x) = \sum_{n\in\bZ}w_ne^{i\xi_n x},\label{eq:fourier}
\end{align}
where
\begin{align*}
    \xi_n = \frac{2\pi n}{L},\quad 
    w_n:= \dfrac{1}{L} \int^{L}_{0} w(x) e^{-i\xi_n x}dx. 
\end{align*}
For $w\in H_{\mathrm{per}}^2(I)$, the Fourier series expansion of $\cL w$ is given by
\begin{align}
     (\cL w)(x)= \sum_{n\in\bZ}\left\{(-d + k\hat{K}(\xi_n))\xi_n^2 - c\right\}w_n e^{i\xi_nx}  
     =: 
     \sum_{n\in\bZ} \lambda_n(k;d,c)w_n e^{i\xi_n x},\label{def:lambda}
\end{align}
where $\hat K$ denotes the Fourier transform of $K$, defined by
$$
\hat{K}(\xi) := \int_{\bR}K(x)e^{-i\xi x}dx.
$$
Since $K$ is even, $\hat K$ is a real-valued even function.
Hence, all $\lambda_n(k;d,c)$ are real.

We next identify possible bifurcation points from $u=1$.
First, the spectrum of $\cL$ is given by
\beaa
\sigma(\mathcal{L})=\{\,\lambda_n(k;d,c)\,\}_{n\ge 0},
\eeaa
and it is immediate that every element of the spectrum is an eigenvalue.
Therefore, a bifurcation of stationary solutions from $u = 1$ may occur when one of the eigenvalues $\lambda_n(k; d, c)$ crosses zero.

\begin{assumption}\label{assumption:lam}
For fixed constants $d,\ c > 0$, there exist $k^* > 0$ and a positive integer $N$ such that, for every positive integer $n$, 
    \begin{align*}
        \lambda_n^* := \lambda_n(k^*;d,c) 
        \begin{cases}
        =0 & \mathrm{if} \ n=N,\\
        \neq 0 & \mathrm{if} \ n\neq N.
        \end{cases} 
    \end{align*}
\end{assumption}
Since $d,\ c$ are fixed, we write the eigenvalues simply as $\lambda_n(k)$.


In this paper, a Turing bifurcation refers to a bifurcation at the first loss of stability of the homogeneous stationary solution $u=1$.
We therefore distinguish the case
\begin{align}\label{condi:first}
    \lambda_n^* := \lambda_n(k^*;d,c) 
    \begin{cases}
    =0 & \mathrm{if} \ n=N,\\
    < 0 & \mathrm{if} \ n\neq N.
    \end{cases} 
\end{align}
Under condition \eqref{condi:first}, it holds that $\hat K(\xi_N)>0$. 
Hence, $\lambda_N(k)$ crosses zero from negative to positive as $k$ increases through $k^*$.
Thus, $k^*$ is a candidate point for a Turing bifurcation.
The existence result below is proved under the more general Assumption~\ref{assumption:lam}. 
Unlike condition~\eqref{condi:first}, this assumption does not require the other eigenvalues to be negative. 
Examples satisfying this assumption are given in Section~\ref{sec:num}.


We first establish the existence of nonhomogeneous stationary solutions arising through bifurcation.
The stationary problem is invariant under spatial translations. 
To make the critical eigenvalue simple, we restrict the problem to the following symmetric function spaces:
\be
\begin{aligned}
    \mathcal{X}
    &:=
    \{\,u \in H_{\mathrm{per}}^2(I)
    \mid u(L-x)=u(x),\ x\in I\,\}, \\
    \mathcal{Y}
    &:=
    \{\,u \in L_{\mathrm{per}}^2(I)
    \mid u(L-x)=u(x),\ x\in I\,\}.
\end{aligned}
\label{eq:function-spaces}
\ee
Let $\tilde{\cL}:\cX \to \cY$ denote the restriction of $\cL$ to $\cX$: 
\beaa
    \tilde{\cL} w = \cL w, \quad w \in \cX.
\eeaa
Note that the spectrum of $\tilde{\cL}$ is given by $\{\lambda_n(k)\}_{n\ge 0}$.
Under Assumption~\ref{assumption:lam}, the kernel of $\tilde{\cL}(k^*)$ is spanned by $\cos(\xi_Nx)$.


By applying the Crandall--Rabinowitz theorem \cite[Theorem 1.7]{CR71}, we obtain the following result on the bifurcation of stationary solutions from $u = 1$:
\begin{proposition}\label{prop:cr}
    Under Assumption~\ref{assumption:lam}, $(k, w) = (k^*, 0)$ is a bifurcation point of $\cH(k, w; d, c) = 0$.
    Moreover, there exists $s^* > 0$ such that in a neighborhood of this point, there exists a set of nonhomogeneous solutions 
    \begin{gather*}
    \{\,(k, w) = (k(s), w(s)) \mid 0 < |s| < s^*\,\},\qquad
    k(0) = k^*, \ w(s) = s\cos(\xi_N x) + o(s), \ s \to 0.
    \end{gather*}
\end{proposition}

\medskip

Proposition~\ref{prop:cr} establishes a local branch of nonhomogeneous stationary solutions bifurcating from $u=1$ at $k=k^*$. 
The result requires Assumption~\ref{assumption:lam}, but not the stronger condition~\eqref{condi:first}. 
Thus, the bifurcation need not coincide with the first loss of stability of $u=1$.
The proof is based on the Crandall--Rabinowitz theorem and is given in Appendix~\ref{sec:A}.

\bigskip
\subsection{Center manifold reduction}\label{sec:cmt}
Proposition~\ref{prop:cr} establishes the existence of a bifurcation of stationary solutions from $u = 1$.
However, it does not determine whether the bifurcation is supercritical or subcritical. 
In this subsection, we derive an amplitude equation using center manifold reduction. 
The sign of its cubic coefficient determines the criticality of the bifurcation.

Before carrying out the reduction, we state local and global existence results for equation~\eqref{eq:main}.
For a general reaction term, the following local existence result holds:
\begin{proposition}\label{prop:local_ex}
    Suppose that $d,k>0$, $K\in W^{1,1}(\bR)$ and $f:\bR\to\bR$ is a locally Lipschitz function satisfying $f(0)=0$.
    Let $\phi \in C_{\mathrm{per}}(I)$.
    Then, there exists $T>0$ such that Equation~\eqref{eq:main}  with $u(0)=\phi$ admits a unique classical solution 
    \beaa
    u\in C([0,T];C_{\mathrm{per}}(I)) \cap C^1((0,T);C_{\mathrm{per}}(I)) \cap C((0,T];C^2_{\mathrm{per}}(I)).
    \eeaa
\end{proposition}

\medskip
Moreover, when the reaction term is logistic, we obtain the following global-in-time existence result:
\begin{proposition}\label{prop:global_ex}
    Suppose that $d,k>0$ and $K\in W^{2,1}(\bR)$.
    Let $\phi \in C_{\mathrm{per}}(I)$ be a nonnegative function.
    Then, Equation~\eqref{eq:main}, with the reaction term \eqref{logistic} and $u(0)=\phi$, admits a unique classical solution
    \beaa
    u\in C([0,+\infty);C_{\mathrm{per}}(I)) \cap C^1((0,+\infty);C_{\mathrm{per}}(I)) \cap C((0,+\infty);C^2_{\mathrm{per}}(I)).
    \eeaa
    Moreover, $u(t)$ is a positive function for $t>0$ if $\phi \not\equiv 0$.
\end{proposition}

\medskip
For Proposition~\ref{prop:cr}, we assumed $K \in L^1(\bR)$ to establish the existence of stationary solutions.
Note, however, that Propositions~\ref{prop:local_ex} and \ref{prop:global_ex} require additional assumptions on the derivatives of $K$.
These assumptions are needed to estimate the nonlocal advection term without imposing unnecessary differentiability assumptions on the solution, and are therefore technical.


For Proposition~\ref{prop:local_ex}, we first define a mild solution and establish its local existence by a fixed-point theorem.
We then show that the mild solution is in fact a classical solution by standard regularity arguments.
For Proposition~\ref{prop:global_ex}, we derive a priori bounds by the energy method and use the resulting boundedness of the solution to extend the local solution globally in time.
Since both proofs are based on standard arguments, they are given in  Appendix B.

We now establish the existence of a local center manifold and derive the reduced dynamics near the bifurcation point.
\begin{theorem}\label{thm:reduction}
    Consider equation \eqref{eq:main_expand} in the symmetric function space $\cX$ defined in \eqref{eq:function-spaces}.
    Suppose that Assumption~\ref{assumption:lam} holds.
    Then there exists a locally invariant center manifold $W^c_{\mathrm{loc}}$ in a neighborhood of the bifurcation point $(k, w) = (k^*, 0)$.
    Moreover, for sufficiently small $\delta > 0$, in the regime
    $$
        |w_N|=O(\delta),\qquad |k-k^*|= O(\delta^2),
    $$
    the reduced equation on $W^c_{\mathrm{loc}}$ is given by
    \begin{align}
        \dot{w}_N
        =\lambda_N(k) w_N  + \mu w^3_N + O(\delta^4), \label{eq:redu}
    \end{align}
    where the coefficient of the cubic term is given by
    \begin{align*}
        \mu 
        = -2(d\xi^2_N-c)
        - \frac{2d\xi^2_N + c}{2\lambda_{2N}(k^*)}\left(\lambda_{2N}(k^*) + 2d\xi^2_N - 5c\right).  
    \end{align*}
\end{theorem}
\begin{proof}
    Suppose that Assumption~\ref{assumption:lam} holds. 
    To focus on the bifurcation of the $N$-th mode at $k = k^*$, we set
    $$
        \kappa := k - k^*,\quad 
        \mathcal{L}^* := \tilde{\cL}(k^*),\quad 
        \lambda_n^* :=\lambda_n(k^*)\ (n\ge 0)
    $$
    and rewrite equation~\eqref{eq:main_expand} as
    \begin{align*}
        \dot{w} = \cL^* w+\cR(w,\kappa),
    \end{align*}
    where
    $$
        \cR(w,\kappa)
        :=
        \bigl(\cL(k^*+\kappa)-\mathcal L^*\bigr)w
        +\cN(w,k^*+\kappa).
    $$
    Since the homogeneous stationary solution $w = 0$ is a solution of the equation for every $\kappa$, it follows that
    $$
        \mathcal R(0,\kappa)=0.
    $$
    In what follows, we apply the parameter-dependent center manifold theorem in \cite[Theorem~3.3]{HI11}.
    This is equivalent to applying the center manifold theorem \cite[Theorem~2.9]{HI11} to the extended system
    \begin{align*}
        \begin{cases}
            \dot{w}
            &=
            \mathcal L^*w+\mathcal R(w,\kappa), \\
            \dot{\kappa}
            &=0.
        \end{cases}
    \end{align*}
    Therefore, to establish the existence of a center manifold, it suffices to verify the following properties:

    \begin{enumerate}
        \item[(C1)]\ $ \cL^*:\cX\to\cY$ is a bounded linear operator.
        Moreover, $\cR(w, \kappa)$ is of class $C^{2}$ in a neighborhood of $0 \in \cX, \kappa = 0$, and satisfies
        $$
            \cR(0,\kappa)=0,
            \quad
            D_w\cR(0,0)=0.
        $$
        Here, $D_w\cR$ denotes the Fr\'echet derivative of $\cR$ with respect to $w$.

        
        \item[(C2)]\ Let
        $
            \sigma_{+}(\cL^*), \
            \sigma_{0}(\cL^*),\
            \sigma_{-}(\cL^*)
        $
        denote the subsets of the spectrum $\sigma(\cL^*)$ consisting of elements with positive, zero, and negative real parts, respectively.
        Then neither $\sigma_{+}(\cL^*)$ nor $\sigma_{-}(\cL^*)$ contains a subsequence converging to the imaginary axis.
        Moreover, $\sigma_0(\cL^*)$ consists of finitely many eigenvalues, counted with multiplicity.

                
        \item[(C3)] For any $\omega>0$, the following estimate holds:
        \begin{align}\label{C3}
            \|(i\omega \mathbb{I} - \cL^*)^{-1}\|_{B(\cY)} 
            \le 
            \frac{1}{\omega}.
        \end{align}
        Here, $\mathbb{I}:\cY\to \cY$ denotes the identity operator, and
        $B(\cY) $ denotes the space of bounded linear operators from $ \cY $ to $ \cY $.
    \end{enumerate}


    We now verify conditions (C1)--(C3).

    
    \textbf{Verification of (C1).} 
    Since the linearity of $\cL^*:\cX\to \cY$ is clear, it remains to verify its boundedness.
    For $w\in \cX$, Young's inequality gives
    \beaa
        \|\cL^* w\|_{L^2(I)} &\le& d \|w_{xx}\|_{L^2(I)} + c \|w\|_{L^2(I)} + k^* \|K\|_{L^1(\bR)}  \|w_{xx}\|_{L^2(I)} \\
        &\le& (d + c + k^* \|K\|_{L^1(\bR)})  \|w\|_{H^2(I)}.
    \eeaa
    Hence, $\cL^*$ is bounded.

    
    It is immediate that $\cR$ is of class $C^2$ in a neighborhood of $(w,\kappa)=(0,0)$ and that $\cR(0,\kappa) = 0$.
    Moreover, for each $h\in \cX$, the Fr\'echet derivative with respect to $w$ is given by
    \beaa
         D_w\cR(w,0)[h] = -2 c w h - k^* [h (K*w)_{x}]_{x} - k^* [w(K*h)_{x}]_{x}
    \eeaa
    and therefore $D_w\mathcal R(0,0)=0$.

    \medskip
    \textbf{Verification of (C2).} 
    We have $\sigma(\cL^*) = \{\,\lambda_n^*\,\}_{n \ge 0}$.
    Since $\lambda_n^*\to -\infty\ (n \to \infty)$, neither $\sigma_+(\cL^*)$ nor $\sigma_-(\cL^*)$ contains a subsequence converging to the imaginary axis.
    Moreover, Assumption~\ref{assumption:lam} gives $\sigma_{0}(\cL^*) = \{\,0\,\}$.
    Since $\ker\cL^* = \mathrm{span}\{\,\cos(\xi_N x)\,\}$, the eigenvalue $0$ is simple.

    \medskip
    \textbf{Verification of (C3).}
    Fix any $\omega>0$.
    For $z \in \cY$, consider
    \be\label{eq:cm-c}
    (i\omega \mathbb{I} - \cL^*)w = z.
    \ee
    Since $i\omega\not\in\sigma(\cL^*)$, there exists a unique $w \in \cX$ satisfying \eqref{eq:cm-c}.
    In particular, if $z\neq 0$, then $w\neq 0$.

    We now prove \eqref{C3}.
    Taking the inner product of both sides of \eqref{eq:cm-c} with $w$ and then taking the imaginary part, we obtain
    \beaa
        \omega\|w\|^2_{L^2(I)} 
        = \mathrm{Im}\int_{I}\, z(x)\overline{w}(x)\,dx.
    \eeaa
    Applying the Cauchy--Schwarz inequality to the right-hand side gives
    \begin{align*}
        \omega\|w\|^2_{L^2(I)} 
        \le \|w\|_{L^2(I)}\|z\|_{L^2(I)}.
    \end{align*}
    Therefore, if $z\neq 0$, then $\|w\|_{L^2(I)}\neq 0$, and dividing both sides by $\|w\|_{L^2(I)}$ yields
    \begin{align*}
         \omega\|w\|_{L^2(I)} \le \|z\|_{L^2(I)}.
    \end{align*}
    Hence, \eqref{C3} follows.


    Therefore, \cite[Theorem~2.9]{HI11} applies, and there exists a locally invariant center manifold $W^c_{\mathrm{loc}}$.


    It remains to compute the reduced equation. 
    The equation for the critical Fourier coefficient $w_N$ is
    \begin{align*}
        &\dot{w}_N
        =\lambda_N(k) w_N + \sum_{m\in\bZ}M_{N,m}(k)w_{m}w_{N-m}, \quad
        &M_{n,m}(k) := k\xi_n\xi_{n-m}\hat{K}(\xi_{n-m})-c.
    \end{align*}
    Under Assumption~\ref{assumption:lam}, we derive the reduced equation using the representation formula for the center manifold \cite[Corollary~2.12]{HI11} and a Taylor expansion.

    
    On the center manifold, the noncritical Fourier coefficients are of second order or higher in $w_N$.
    Therefore, for sufficiently small $\delta>0$, under the scaling $w_{N} = O(\delta), \ |k-k^*| = O(\delta^2) $, we have
    \beaa
        w_n = O(\delta^2)\quad (n\neq \pm N).
    \eeaa
    Using
    \beaa
        \lambda^{*}_N = - d \xi^2_N + k^* \xi^2_N \hat{K}(\xi_N) - c = 0
    \eeaa
    and
    \beaa
        \lambda^{*}_{2N} 
        = 
        \bigl(- d  + k^* \hat{K}(\xi_{2N})\bigr)\xi^2_{2N} - c,
    \eeaa
    we obtain
    \begin{align}
        \dot{w}_N
        =&\lambda_N(k) w_N + \sum_{m\in\bZ}M^*_{N,m}w_{m}w_{N-m} +O(\delta^4)\notag \\
        =&\lambda_N(k) w_N + (M^*_{N,N} + M^*_{N,0})w_0 w_N + (M^*_{N,-N} + M^*_{N,2N})w_{2N}w_{N}+ O(\delta^4) \notag\\
        =&\lambda_N(k) w_N + (d \xi^2_N -c)w_0 w_N + \dfrac{1}{2} ( \lambda^*_{2N} +2d\xi^2_N - 5c)w_{2N}w_{N}+ O(\delta^4) \label{eq:uN}
    \end{align}
    where
    $$
        M_{n,m}^* := k^*\xi_n\xi_{n-m}\hat{K}(\xi_{n-m})-c.
    $$
    We also note that $w_{-n} = w_n$ because we restrict the problem to the space of even functions.


    Next, the center manifold approximation gives
    \begin{align}
        w_0 = a_0 w_N^2 + O(\delta^3),\quad
        w_{2N} = a_{2N}w_N^2 + O(\delta^3). \label{eq:u0}
    \end{align}
    It therefore remains to determine $a_0$ and $a_{2N}$.
    Substituting \eqref{eq:u0} into the invariance equation for the center manifold and comparing the coefficients of $w_N^2$, we obtain
    \begin{align*}
        \lambda_0(k^*)a_0+M^*_{0,N}+M^*_{0,-N}=0,\quad
        \lambda_{2N}(k^*)a_{2N}+M^*_{2N,N}=0.
    \end{align*}
    Since $\lambda_0(k^*)=-c$ and $M^*_{0,N}+M^*_{0,-N}=-2c$, we obtain $a_0=-2$.
    Moreover, Assumption~\ref{assumption:lam} gives $\lambda_{2N}(k^*)\neq0$, and $M^*_{2N,N}=2d\xi_N^2+c$.
    Hence,
    $$
        a_{2N}=-\frac{2d\xi_N^2+c}{\lambda_{2N}(k^*)}.
    $$
    Substituting these expressions into \eqref{eq:uN} gives the reduced equation \eqref{eq:redu}.
    This completes the proof.
\end{proof}


First, we note a property of the center manifold.
The center manifold is locally invariant, and for suitable initial data, the behavior of solutions can be described by the reduced equation.
In particular, under condition \eqref{condi:first} for the Turing bifurcation, the spectrum of $\cL^*$ satisfies $\sigma_+(\cL^*)=\emptyset$.
Therefore, \cite[Theorem~3.22]{HI11} gives the following result:
\begin{corollary}\label{cor:attract}
    If \eqref{condi:first} holds, then $W^c_{\mathrm{loc}}$ is locally attracting.
\end{corollary}

\medskip

Corollary~\ref{cor:attract} shows that the center manifold attracts nearby trajectories in the symmetric function space. 
The reduced equation therefore describes the asymptotic behavior of nearby trajectories that remain in a sufficiently small neighborhood.

For $\mu\neq0$, equation~\eqref{eq:redu} describes a pitchfork bifurcation from $u=1$. 
The bifurcation is supercritical if $\mu<0$ and subcritical if $\mu>0$.
Under condition~\eqref{condi:first}, the supercritical branch consists of stable nonhomogeneous stationary solutions near $k=k^*$.
In the subcritical case, the bifurcating solutions are unstable and occur where the homogeneous solution is stable. 
These stability statements refer to perturbations in the symmetric function space.
A subcritical Turing bifurcation may lead to multistability. 
This occurs if the bifurcating branch becomes stable through a saddle-node bifurcation while the homogeneous solution remains stable. 
In Section~\ref{sec:deformation}, we examine how the kernel affects the condition for subcriticality.


\subsection{Conditions for subcritical bifurcation}\label{sec:deformation}
In this subsection, we express the subcriticality condition $\mu>0$ in terms of the interaction kernel. 
Throughout, we assume that the bifurcation point satisfies condition~\eqref{condi:first}. 
The criterion is written in terms of the critical mode $N$ and the ratio $\hat{K}(\xi_{2N})/\hat{K}(\xi_N)$.

\begin{figure}[tb!]
  \centering
    \includegraphics[width=1.0\linewidth]{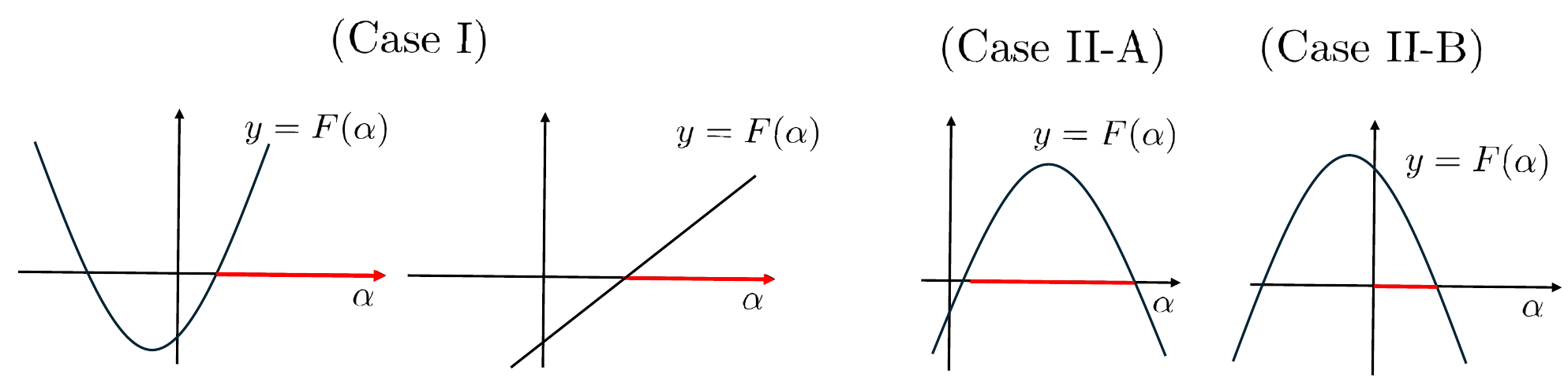}
  \caption{\small{
  The red intervals indicate the values of $\alpha =d\xi^2_N>0$ for which $F(\alpha)>0$.
  }}
  \label{fig:case1_2}
\end{figure}

First, note that \eqref{condi:first} implies $\lambda^*_{2N} < 0$.
Then the condition $\mu > 0$ for a subcritical bifurcation can be written as
\beaa
    0 < 3\lambda^*_{2N}(2d\xi^2_N-c)+(2d\xi^2_N+c)(2d\xi^2_N-5c).
\eeaa
Using $\lambda_N^* = 0$, we also have
\begin{align*}
    k^* 
    = \left(d +\frac{c}{\xi_N^2}\right) \hat{K}(\xi_N)^{-1} 
    > 0,
\end{align*}
and hence the inequality can be rewritten as
\begin{align}
    \frac{10(d\xi^2_N)^2+d\xi^2_Nc+c{}^2}{6} 
    < (2d\xi^2_N-c)(d\xi^2_N+c)\frac{\hat{K}(\xi_{2N})}{\hat{K}(\xi_{N})}.\label{ineq:subcri}
\end{align}
This criterion makes the dependence on the Fourier transform of the interaction kernel explicit.
Furthermore, setting
\begin{align*}
    \alpha := d\xi^2_N > 0,\quad
    R := \frac{\hat{K}(\xi_{2N})}{\hat{K}(\xi_{N})},
\end{align*}
we can rewrite \eqref{ineq:subcri} as
\begin{align}
    0 < 2(6R - 5)\alpha^2 + (6R - 1)c\alpha - c^2(6R + 1) 
    =: F(\alpha). \label{F(x)>0}
\end{align}
We now determine when $F(\alpha)>0$ for $\alpha>0$.
The analysis is divided into two cases according to the sign of $6R-5$:

\medskip
\textbf{Case I.} 
    Suppose that $6R - 5 \ge 0$.
    Since $F(0) = -c^2(6R + 1) < 0$, equation $F(\alpha) = 0$ has exactly one solution $\alpha^* > 0$.
    Moreover, inequality \eqref{F(x)>0} holds for $\alpha>\alpha^*$.
    Therefore, as shown in the Case~I panel of Fig.~\ref{fig:case1_2}, if $6R - 5 \ge 0$ and $\alpha$ is sufficiently large, the bifurcation is subcritical.

\textbf{Case II.} 
    Suppose that $6R - 5 < 0$.
    Let $D(R)$ denote the discriminant of the quadratic equation $F(\alpha) = 0$.
    Then
    $$
        D(R) = 3c^2 (108R^2 - 68R - 13).
    $$
    For $F(\alpha)>0$ to hold, it is necessary that $D(R)>0$.
    The solutions of $D(R)=0$ are given by
    $$
        R_{{\pm}} = \frac{17 \pm 8\sqrt{10}}{54},
    $$
    and hence $D(R)>0$ when
    \begin{align}\label{condi:R}
        R < R_{-}, \quad R_{+} < R< \frac{5}{6}.
    \end{align}
    
    
    For $F$ to be positive at some $\alpha>0$, one of the following two cases must hold:
    \begin{description}
        \item[Case II-A] 
        Suppose that $D(R) > 0$ and $F(0) < 0$.
        We also require the critical point $\alpha_c$ of $F$ to satisfy $\alpha_c>0$.
        That is, \eqref{condi:R} and
        \beaa
            F(0) = -c^2(6R + 1) < 0 
        \eeaa
        and
        \beaa
            \alpha_c = -\frac{c(6R - 1)}{4(6R - 5)} > 0
        \eeaa
        hold.
        Since $R_+>\frac{1}{6}$, these conditions reduce to $R_{+} < R<\frac{5}{6}$.


        \item[Case II-B] 
        Suppose that $D(R)>0$ and $F(0)>0$. 
        Since $F(0)=-c^2(6R+1)$, these conditions reduce to $R<-1/6$.

    \end{description}
As shown in the Case~II-A and Case~II-B panels of Fig.~\ref{fig:case1_2}, the range of $\alpha$ satisfying inequality \eqref{F(x)>0} is bounded.

\medskip

The above results allow us to interpret the condition for a subcritical Turing bifurcation in terms of the relation between the Fourier transform of the interaction kernel and the critical mode $N$.
Moreover, condition \eqref{ineq:subcri} not only determines the criticality for a given interaction kernel but also provides a guideline for constructing interaction kernels that produce a subcritical Turing bifurcation.


The authors of \cite{BT21} performed a weakly nonlinear analysis of a model without a reaction term. 
They showed that the bifurcation direction depends on the interaction kernel.
They also constructed specific interaction kernels satisfying a sign condition for subcriticality.

In Section~\ref{subsec:rescale}, we consider a parameterized family of interaction kernels and show that the above condition for a subcritical bifurcation is satisfied for a particular class of kernels.
This result corresponds to Case~I.
In Section~\ref{sec:num}, we construct specific interaction kernels satisfying the derived conditions and numerically confirm the occurrence of subcritical bifurcations and multistability between homogeneous and nonhomogeneous stationary solutions.
Cases~II-A and~II-B are investigated in Sections~\ref{subsec:caseIIa} and~\ref{subsec:caseIIb}, respectively.


\subsection{Subcritical bifurcation induced by interaction kernel rescaling}\label{subsec:rescale}
In this subsection, we examine how interaction kernel rescaling affects the criticality of the Turing bifurcation. 
Under the assumptions below, sufficiently localized kernels satisfy the subcriticality condition in Case~I of Section~\ref{sec:deformation}.

\begin{theorem}\label{thm:deformation}
    Let $J\in L^{1}(\bR)$ be an even function satisfying
    \beaa
        \hat{J}(0)>0,\quad \hat{J}(0)>\hat{J}(\xi) \quad (\xi\neq 0).
    \eeaa
    For $\beta>0$, define
    \be\label{ker:rescale}
        K(x) = \frac{1}{\beta} J\left(\frac{x}{\beta}\right).
    \ee
    
    Suppose that $\{\beta_j\}_{j\in\bN}$ is a sequence of positive numbers with $\beta_j\to0$.
    For each $j$, assume that condition \eqref{condi:first} holds for the kernel \eqref{ker:rescale} with $\beta=\beta_j$.
    Let $k_j^*>0$ and $N_j\in\bN$ denote the corresponding bifurcation point and critical mode.
    Then the cubic coefficient $\mu_j$ in the reduced equation of Theorem~\ref{thm:reduction} satisfies
    $\mu_j > 0$
    for all sufficiently large $j$.
\end{theorem}

\medskip
Under the stated assumptions, localization of the interaction leads to a subcritical Turing bifurcation along the sequence $\beta_j\to0$.
For example, the Fourier-transform conditions hold if $J$ is nonnegative, continuous, integrable, even, and not identically zero.


We first prepare some notation for the proof of Theorem~\ref{thm:deformation}.
By the definition of $K$, its Fourier transform satisfies $\hat{K}(\xi)=\hat{J}(\beta\xi)$.
Fix $d,c>0$.
For each $n\in\bN$, define
\begin{align*}
    &k^+(n;\beta) :=
    \begin{cases}
        \dfrac{1}{\hat{J}(\beta\xi_n)}\left(d + \dfrac{c}{\xi_n^2}\right), &(\hat{J}(\beta\xi_n) >0),\\
         +\infty,  &(\hat{J}(\beta\xi_n) \le 0).
    \end{cases}
\end{align*}
Here, $k^+(n;\beta)$ is the instability threshold for the $n$-th mode. 
The value $+\infty$ indicates that this mode remains stable for all \(k>0\).


At $\beta=\beta_j$, the critical mode is given by
\beaa
    N_j :=\underset{n \in \mathbb{N}}{\operatorname{argmin}} \ k^+(n;\beta_j).
\eeaa
This minimizer is unique by condition~\eqref{condi:first}.

For each fixed $\xi_n$, we also define the limit of $k^+(n;\beta)$ as $\beta\to0$ by
\beaa
    k_0(n): =\dfrac{1}{\hat{J}(0)} \left(d + \frac{c}{\xi_n^2}\right).
\eeaa
By the assumption of Theorem~\ref{thm:deformation}, for sufficiently small $\beta$ we have $\hat{J}(\beta\xi_n)<\hat{J}(0)$.
Hence, $k^+(n;\beta)\ge k_0(n)$.


In the present setting, the condition for subcriticality \eqref{ineq:subcri} can be rewritten as
\begin{align}
    \frac{10(d\xi^2_{N_j})^2+d\xi^2_{N_j}c+c{}^2}{6} 
    < (2d\xi^2_{N_j}-c)(d\xi^2_{N_j}+c)\frac{\hat{J}(\beta_j\xi_{2N_j})}{\hat{J}(\beta_j\xi_{N_j})}. \label{ineq:subcri_b}
\end{align}
We now show that this inequality holds for all sufficiently large $j$.


First, we prove the following lemma on the behavior of $N_j$.
\begin{lemma} \label{lem:sig_inf}
    $\dlim_{j\to \infty}N_j =\infty$.
\end{lemma}
\begin{proof}
    Suppose, for contradiction, that $N_j$ does not tend to infinity.
    Since $N_j\in\bN$, there is a constant subsequence.
    After passing to this subsequence and relabeling, we may assume that
    $N_j=\tilde{N}$ for all $j$, where $\tilde{N}\in\bN$ is fixed.
    Since $\beta_j\to0$, we have
    \begin{equation*}
        \lim_{j\to\infty}\beta_j\xi_{N_j}=0.
    \end{equation*}

    By the minimizing property of $N_j$,
    \beaa
        k^+(\tilde{N};\beta_j) = k^+(N_j;\beta_j) \le k^+(\tilde{N}+1;\beta_j).
    \eeaa
    For each fixed $n\in\bN$, we have
    $k^+(n;\beta_j)\to k_0(n)$ as $j\to\infty$.
    Passing to the limit in the above inequality yields
    \beaa
        k_0(\tilde{N})\le k_0(\tilde{N}+1),
    \eeaa
    which contradicts the strict decrease of $k_0$.
\end{proof}

\medskip

From this lemma, for sufficiently large $j$, the dominant terms on both sides of inequality \eqref{ineq:subcri_b} are those involving $d\xi_{N_j}^2$.
Suppose that there exists some $\tilde J>0$ such that
$$
    \lim_{j\to\infty}\frac{\hat{J}(\beta_j\xi_{2N_j})}{\hat{J}(\beta_j\xi_{N_j})}
    = \tilde J.
$$
Then, as $j\to\infty$,
$$
    \dfrac{\mathrm{(LHS)}}{(d\xi^2_{N_j})^2} \to \frac{5}{3}, \quad \dfrac{\mathrm{(RHS)}}{(d\xi^2_{N_j})^2} \to 2\tilde J.
$$
Therefore, if $\tilde J>5/6$, inequality \eqref{ineq:subcri_b} holds.
To justify this argument, we prepare the following lemma.

\begin{lemma}\label{lem:b_sigma}
    $\dlim_{j\to \infty}\beta_j N_j = 0$.
\end{lemma}
\begin{proof}
    First, we prove that
    \begin{equation}\label{lim}
        \lim_{j\to\infty} k^+(N_j;\beta_j) = \inf_{n\in\bN} k_0(n) = \frac{d}{\hat{J}(0)}.
    \end{equation}
    For each fixed $m\in\bN$, the definitions of $k^+$ and $N_j$ give
    \beaa
        \frac{d}{\hat{J}(0)} \le k_0(N_j) \le k^+(N_j;\beta_j) \le k^+(m;\beta_j).
    \eeaa
    Since $k^+(m;\beta_j)\to k_0(m)$ as $j\to\infty$, we obtain
    \beaa
        \frac{d}{\hat{J}(0)}\le \liminf_{j\to\infty} k^+(N_j;\beta_j) \le \limsup_{j\to\infty} k^+(N_j;\beta_j) \le k_0(m).
    \eeaa
    Letting $m\to\infty$ and using $k_0(m)\to d/\hat{J}(0)$ proves \eqref{lim}.

    Next, set $\xi^*_{N_j}:=\beta_j\xi_{N_j}$.
    By Lemma~\ref{lem:sig_inf} and \eqref{lim},
    \begin{equation}\label{A}
        \hat{J}(\xi^*_{N_j})
        = \frac{d+c/\xi_{N_j}^2}{k^+(N_j;\beta_j)}
        \longrightarrow \hat{J}(0)>0.
    \end{equation}
    Let
    \beaa
        s:=\limsup_{j\to\infty}\xi^*_{N_j}\in[0,+\infty].
    \eeaa
    Suppose, for contradiction, that $s>0$.
    Passing to a subsequence, we may assume that $\xi^*_{N_j}\to s$.

    If $s<+\infty$, continuity of $\hat{J}$ and \eqref{A} imply $\hat{J}(s)=\hat{J}(0)$.
    This contradicts the hypothesis of Theorem~\ref{thm:deformation}.
    If $s=+\infty$, the Riemann--Lebesgue lemma gives $\hat{J}(\xi^*_{N_j})\to0$, contradicting \eqref{A}.
    Hence $s=0$.
    Since $\xi^*_{N_j}\ge0$, it follows that
    \beaa
        \beta_j N_j=\frac{L}{2\pi}\xi^*_{N_j}\to 0.
    \eeaa
    This proves the lemma.
\end{proof}

\medskip
Finally, we prove Theorem~\ref{thm:deformation}.
\begin{proof}[Proof of Theorem~\ref{thm:deformation}]
    By Lemma~\ref{lem:b_sigma},
    $$
        \lim_{j\to\infty}\frac{\hat{J}(\beta_j\xi_{2N_j})}{\hat{J}(\beta_j\xi_{N_j})} 
        = 1 
    $$
    holds.
    By Lemma~\ref{lem:sig_inf}, for sufficiently large $j$, it is enough to compare the coefficients of the highest-order terms in $N$ on both sides of \eqref{ineq:subcri_b}.
    Comparing the coefficients of $\xi^4_{N_j}$ shows that inequality \eqref{ineq:subcri_b} holds.
    Therefore, the cubic coefficient $\mu_j$ in the reduced equation of Theorem~\ref{thm:reduction} satisfies $\mu_j>0$ for all sufficiently large $j$.
\end{proof}

\medskip

We have obtained a sufficient condition for a subcritical Turing bifurcation through a rescaling of the kernel.
As seen from the proof above, the condition in Theorem~\ref{thm:deformation} corresponds to the condition for subcriticality in Case~I of Section~\ref{sec:deformation}.


Theorem~\ref{thm:deformation} assumes that condition~\eqref{condi:first} holds along a sequence $\beta_j\to0$.
Establishing the existence of such a sequence requires a separate argument. For individual values of $\beta$, however, condition~\eqref{condi:first} can be examined numerically. 
In Section~\ref{sec:aggregation}, we numerically investigate subcritical Turing bifurcations as $\beta$ decreases.

\bigskip
\section{Numerical verifications of multistability}\label{sec:num}
In this section, we use numerical bifurcation analysis to investigate multistability associated with subcritical Turing bifurcations. 
The analysis in Section~\ref{sec:sta} determines the bifurcation direction near the bifurcation point. 
However, it does not determine whether the bifurcating branch becomes stable away from that point. 
We therefore trace branches of nonhomogeneous stationary solutions over a wide parameter range. 
Their stability and secondary bifurcations are examined to identify multistability involving stationary and time-periodic solutions.

The computations use a Julia implementation based on \cite{AT72b,AT72a,DKK91a,DKK91b}.
This implementation is used to trace branches of stationary and time-periodic solutions and determine their linear stability.

For the computations, we expand $u$ in a Fourier series as in \eqref{eq:fourier}. 
Substitution into equation~\eqref{eq:main} with the logistic reaction term~\eqref{logistic} gives
\begin{align*}
    \dot{u}_n
    =(-d\xi_n^2+c)u_n(t)+\sum_{m\in \bZ}(k\xi_n\xi_{n-m}\hat{K}(\xi_{n-m})-c)u_m(t)u_{n-m}(t) \quad (n\in\bZ).
\end{align*}
We truncate this system to finitely many Fourier modes.
We retain the Fourier modes with $|n|\le N_{wave}$ and set all coefficients outside this range to zero.
The resulting finite-dimensional system is
\begin{align*}
    \dot{u}_n
    =&(-d\xi_n^2+c)u_n(t)\\
    &+\sum_{|m|\le N_{wave}}(k\xi_n\xi_{n-m}\hat{K}(\xi_{n-m})-c)u_m(t)u_{n-m}(t) \quad (n\in\bZ,\ |n|\le N_{wave}).
\end{align*}
In the following numerical experiments, we set $N_{wave}=32$.

Motivated by Theorem~\ref{thm:deformation}, Section~\ref{sec:aggregation} examines how localization affects the bifurcation structure for three families of attractive interaction kernels.
We then construct settings corresponding to Cases~II-A and II-B in Section~\ref{sec:deformation} and present the bifurcation structures observed in each case.

Unless otherwise stated, the bifurcation diagrams show $k$ on the horizontal axis and the following quantity $A(u)$ on the vertical axis:
\beaa
    \cA(u):=
    \begin{cases}
        \|u\|_{L^2(I)} & (u\ \mathrm{is\ stationary}), \\
        \dmax_{t\in[0,T_p]}\|u(t)\|_{L^2(I)} &\left(u\ \mathrm{ is\ periodic\ with\ period}\ T_p \right).
    \end{cases}
\eeaa


\subsection{Strongly localized attractive interactions}\label{sec:aggregation}


Here, we consider three families of purely attractive interaction kernels of the form \eqref{ker:rescale}.
The corresponding base kernels $J$ satisfy the Fourier-transform conditions in Theorem~\ref{thm:deformation}.
Decreasing $\beta$ makes the interactions more localized. 
We examine how this localization affects the bifurcation structure and the associated multistability.

\subsubsection{Gaussian kernel}\label{subsubsec:gaussian}
\begin{figure}[!tb]
\centering
\includegraphics[width=\linewidth]{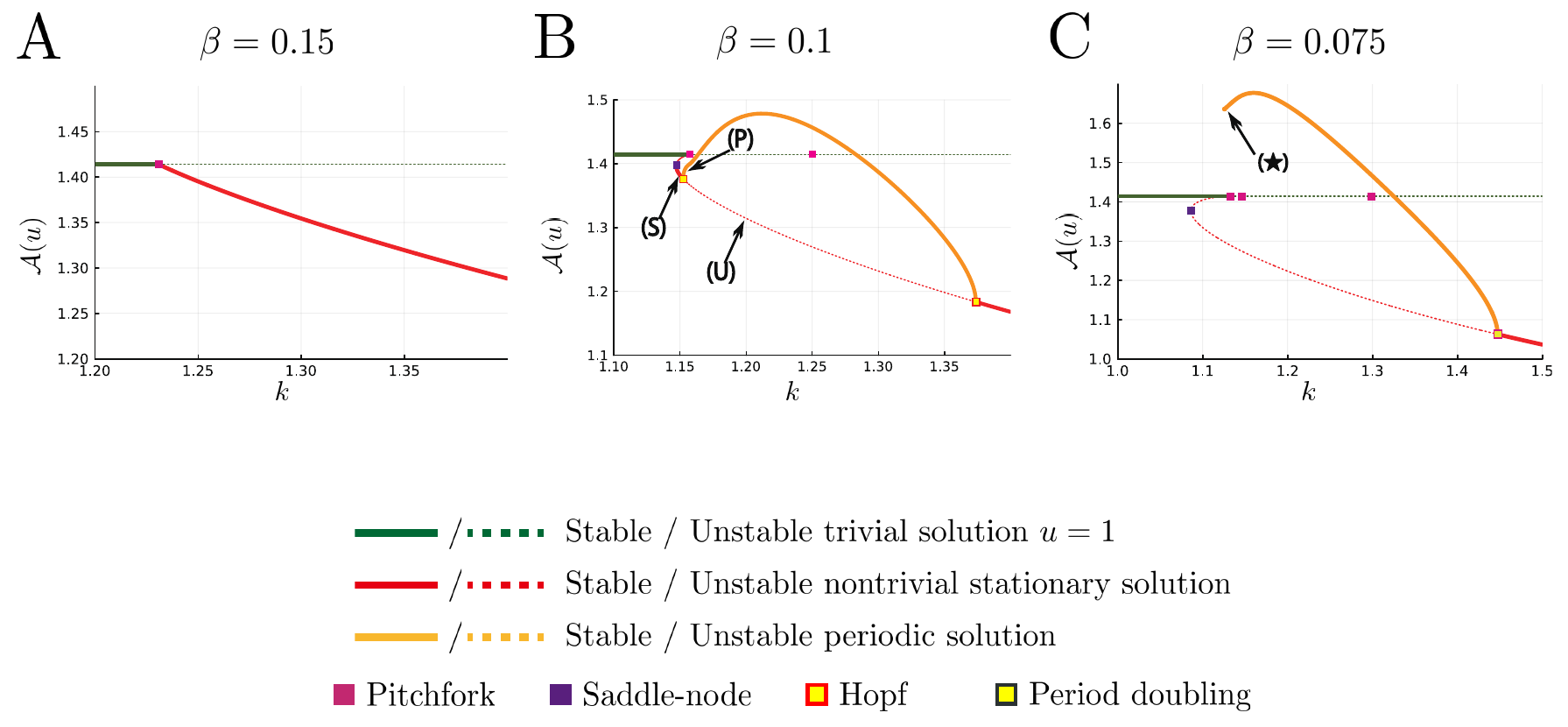}
\caption{\small{
    Bifurcation diagrams for the Gaussian kernel \eqref{Gaussian}.
    We set $d = c = 1.0$ and $L = 2.0$.
    A: Bifurcation diagram for $\beta = 0.15$.
    B: Bifurcation diagram for $\beta = 0.1$.
    C: Bifurcation diagram for $\beta = 0.075$.
    The symbol ($\star$) corresponds to $k=k_{\star} \simeq 1.125645$.}}
\label{fig:gaussian}
\end{figure}

\begin{figure}[!tb]
  \centering
    \includegraphics[width=1.0\linewidth]{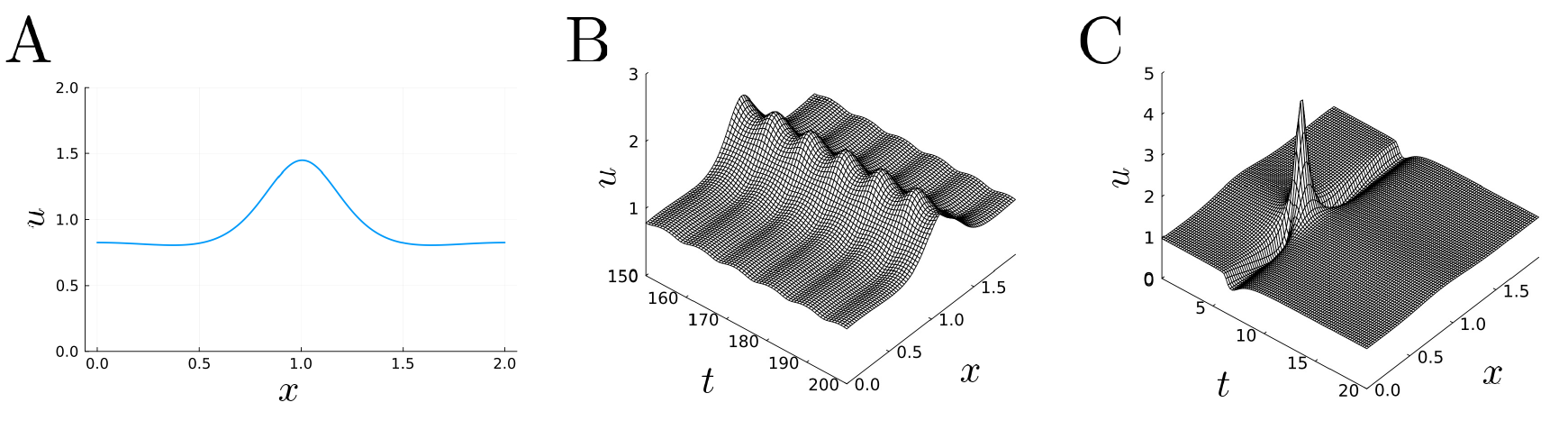}
  \caption{\small{
  Numerical simulations for the Gaussian kernel \eqref{Gaussian}.
  A: A stable stationary solution observed at $k=1.15$ for $\beta=0.1$.
  B: A periodic solution observed at $k=1.153$ for $\beta=0.1$.
  C: Solution dynamics observed at $k=1.1256<k_{\star}$ for $\beta=0.075$.}}
  \label{fig:gaussian_sol}
\end{figure}

We first consider the Gaussian kernel
\be\label{Gaussian}
    K(x)= \frac{1}{\sqrt{2\pi}\beta}\exp\left( -\frac{x^2}{2\beta^2}\right).
\ee
Figure~\ref{fig:gaussian} shows bifurcation diagrams for three values of $\beta$. 
Representative solutions and their time evolution are shown in Figure~\ref{fig:gaussian_sol}.

We choose the parameters so that the first Fourier mode is the critical mode when the homogeneous stationary solution first loses stability.
For $\beta = 0.15$ (Figure~\ref{fig:gaussian}~A), the bifurcation is supercritical.
In contrast, for smaller values of $\beta$, the bifurcation becomes subcritical, and two different bifurcation structures are observed (Figure~\ref{fig:gaussian}~B--C).

\medskip

\textbf{Multistability.}\quad
For $\beta = 0.1$, the nonhomogeneous stationary solutions arising from the subcritical bifurcation become stable through a saddle-node bifurcation (Figure~\ref{fig:gaussian}~B).
As a result, there is a parameter region in which the homogeneous stationary solution and nonhomogeneous stationary solutions are both stable.
One of the stable stationary solutions obtained in this region is shown in Figure~\ref{fig:gaussian_sol}~A.

The stabilized branch of nonhomogeneous stationary solutions then loses stability through a Hopf bifurcation, from which a branch of stable periodic solutions emerges.
Part of this periodic branch lies in the parameter region where the homogeneous stationary solution is stable.
Thus, multistability between the homogeneous stationary solution and time-periodic solutions occurs.
The periodic solution in Figure~\ref{fig:gaussian_sol}~B has a single-peaked spatial profile whose amplitude oscillates in time.
The periodic branch then rejoins the same one-mode branch while remaining stable, and the one-mode branch becomes stable again.

For $\beta = 0.075$, the nonhomogeneous stationary solutions arising from the subcritical bifurcation remain unstable after passing through a saddle-node bifurcation (Figure~\ref{fig:gaussian}~C).
Therefore, multistability between the homogeneous stationary solution and nonhomogeneous stationary solutions does not occur near the bifurcation point.
However, a Hopf bifurcation occurs on the same nonhomogeneous stationary branch near $k=1.45$, away from the bifurcation point, and a branch of stable periodic solutions extends toward smaller values of $k$.
The symbol ($\star$) in Figure~\ref{fig:gaussian}~C corresponds to $k=k_{\star} \simeq 1.125645$, and the stable periodic branch extends at least to $k=k_{\star}$.
The value $k=k_{\star}$ is smaller than the value of $k$ at which the first mode loses stability.
Thus, a parameter region with multistability between the homogeneous stationary solution and stable periodic solutions is also observed for this value of $\beta$.
At $k=1.1256$, the trajectory shown in Figure~\ref{fig:gaussian_sol}~C exhibits transient growth before returning to the homogeneous stationary solution.

\begin{figure}[!tb]
  \centering
    \includegraphics[width=1.0\linewidth]{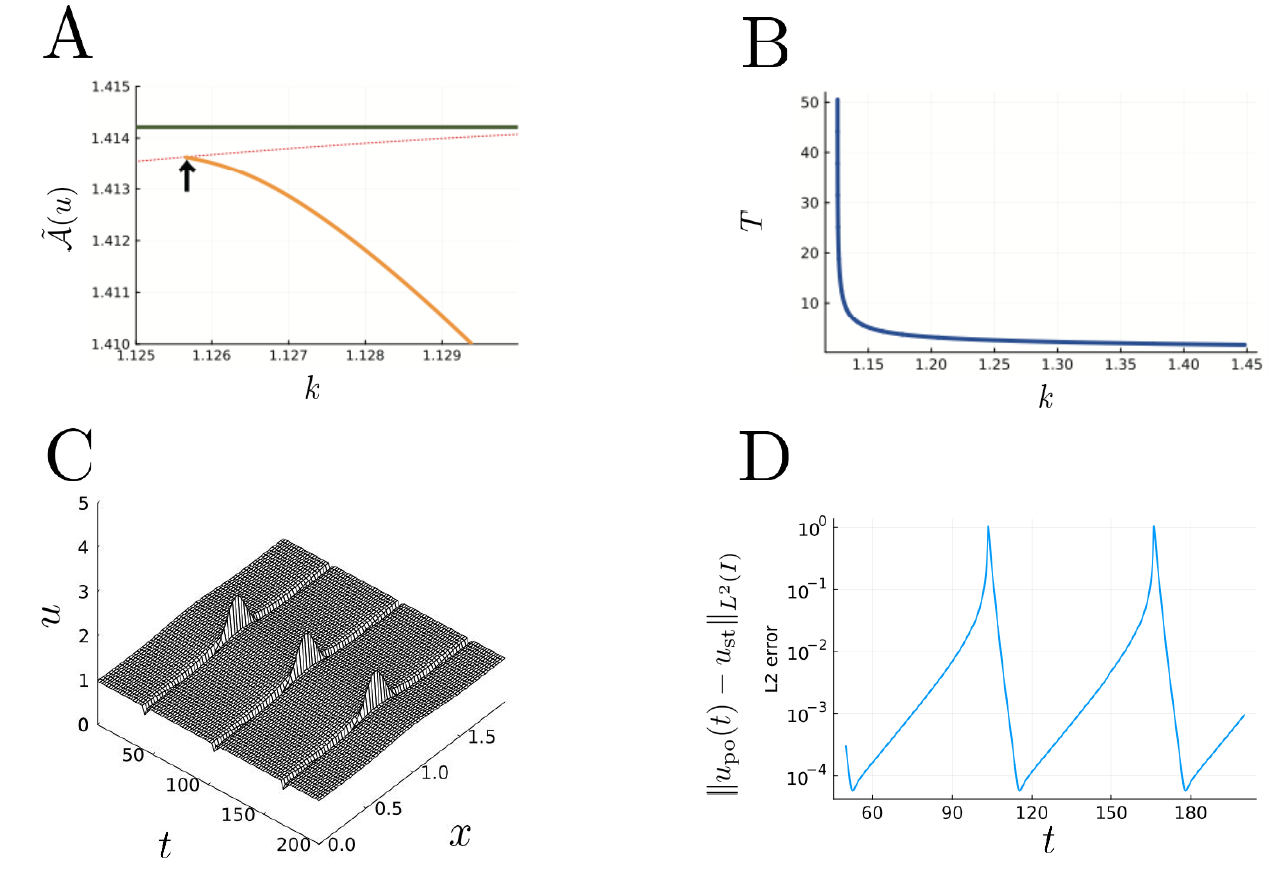}
  \caption{\small{
  A: Bifurcation diagram for the Gaussian kernel \eqref{Gaussian} with $\beta = 0.075$, where the vertical axis is replaced by $\tilde{\cA}(u)$.
  B: Relation between $k$ and the period of the periodic solution for $k>k_{\star}$.
  C: Solution behavior of a periodic solution observed at $k=1.12565$.
  D: $L^2$-distance between the unstable stationary solution at $k=k_{\star}$ and the periodic solution at $k=1.12565$.
  }}
  \label{fig:gaussian_homocli}
\end{figure}

\medskip
\textbf{Numerical evidence for a homoclinic bifurcation.}\quad
For $\beta = 0.075$, we next examine the structure of the periodic solutions near $k=k_{\star}$.
Numerically, we observe that the branch of periodic solutions arising from the Hopf bifurcation near $k=1.45$ enters a parameter region of multistability with the homogeneous stationary solution and does not connect to another Hopf bifurcation point.
Fig.~\ref{fig:gaussian_homocli}~A shows a bifurcation diagram in which the vertical axis is chosen as 
\beaa
    \tilde{\cA}(u):=
    \begin{cases}
        \|u\|_{L^2(I)} & (u\ \mathrm{is\ stationary}), \\
         \|u(t^*)\|_{L^2(I)} &\left(u\ \mathrm{ is\ periodic\ with\ period}\ T_p,\ t^* = \operatorname{argmin}_{t\in [0,T_p]} \left\|u_t(t)\right\|_{L^2(I)} \right).
    \end{cases}
\eeaa
Replacing the vertical axis by $\tilde{\cA}(u)$ makes it possible to see that the branch of periodic solutions approaches the branch of unstable stationary solutions.


Figure~\ref{fig:gaussian_homocli}~B shows the period of the stable periodic solution as a function of the bifurcation parameter $k$.
As $k$ approaches $k_{\star}$, the period increases rapidly.
Figure~\ref{fig:gaussian_homocli}~C shows a space--time plot of a stable periodic solution observed for $k>k_{\star}$ with $k$ close to $k_{\star}$.
The periodic solution stays near a certain state for a long time and then exhibits a rapid large-amplitude excursion.
In Figure~\ref{fig:gaussian_homocli}~D, we plot the $L^2$-distance between this periodic solution and the unstable stationary solution at $k=k_{\star}$.
The periodic solution stays close to the unstable stationary solution for a certain time, then moves away significantly, and later returns close to it again.
These observations suggest that the periodic solutions near $k=k_{\star}$ stay near the unstable stationary solution for a long time.
The increasing period and repeated returns near the unstable stationary solution suggest that the periodic branch approaches a homoclinic bifurcation near $k=k_\star$.
This interpretation is based on numerical evidence.

\medskip

\subsubsection{Top-hat kernel}
Next, we consider the following interaction kernel, which also represents attractive interactions:
\begin{align}\label{top-hat}
    K(x)=
    \left\{
        \begin{aligned}
        &\frac{1}{\beta}, \quad x \in \left[-\frac{\beta}{2}, \frac{\beta}{2}\right], \\
        &0,\quad x \notin \left[-\frac{\beta}{2}, \frac{\beta}{2}\right].
        \end{aligned}
    \right.
\end{align}
We note that this kernel does not satisfy the regularity assumption
$K\in W^{1,1}(\mathbb R)$ imposed in Proposition~\ref{prop:local_ex}.
The velocity induced by the top-hat kernel depends on the density difference between $x-\beta/2$ and $x+\beta/2$, directing movement toward the higher-density side.
Indeed, the term corresponding to the velocity field in the nonlocal advection can be written as
$$
    (K*u)_x(t, x) 
    = (K*u_x)(t, x) 
    = \frac{1}{\beta}\int_{x - \beta/2}^{x + \beta/2}u_{x}(t, y)dy 
    = \frac{u(t, x + \beta/2) - u(t, x - \beta/2)}{\beta}.
$$

\begin{figure}[!tb]
  \centering
    \includegraphics[width=\linewidth]{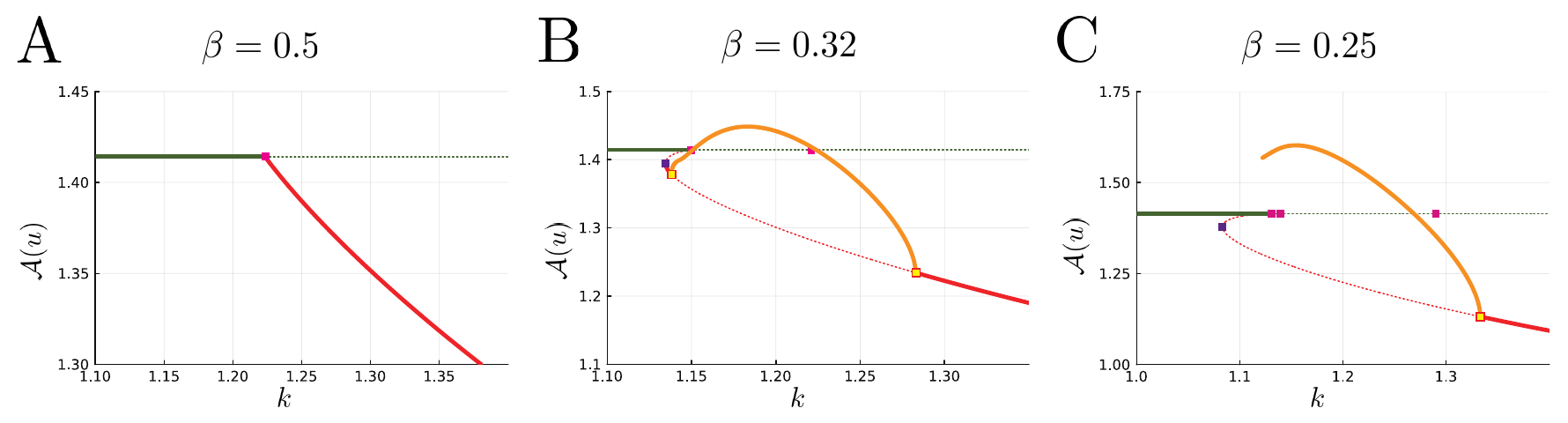}
  \caption{\small{Bifurcation diagrams for the top-hat kernel \eqref{top-hat}.
      See Figure~\ref{fig:gaussian} for the correspondence between branch colors and solution types.
      We set $d = c = 1.0$ and $L = 2.0$.
      A: Bifurcation diagram for $\beta = 0.5$.
      B: Bifurcation diagram for $\beta = 0.32$.
      C: Bifurcation diagram for $\beta = 0.25$.
  }}
  \label{fig:indicator}
\end{figure}

Figure~\ref{fig:indicator} shows the bifurcation diagrams for three different values of $\beta$.
As in the Gaussian kernel considered in Section~\ref{subsubsec:gaussian}, we choose the parameters so that the first Fourier mode is the first mode to bifurcate.
For $\beta=0.5$, a supercritical bifurcation is observed (Figure~\ref{fig:indicator}~A).
For smaller values of $\beta$, two characteristic bifurcation structures associated with subcritical bifurcations are observed (Figure~\ref{fig:indicator}~B--C).
Both structures are qualitatively similar to those in Figure~\ref{fig:gaussian}~B--C.
For the case $\beta=0.25$ shown in Figure~\ref{fig:indicator}~C, the stable periodic branch extends toward smaller values of $k$ without reconnecting to another Hopf bifurcation point in the computed range.
Near the end of the computed branch, the periodic behaviour suggests a homoclinic bifurcation, as in the Gaussian case.

\medskip

\subsubsection{Exponential kernel}
Finally, we consider
\be\label{exp-ker}
    K(x) = \frac{1}{2\beta}\exp\left(-\frac{|x|}{\beta}\right).
\ee
For this kernel, setting $v(t,x):= (K*u)(t,x)$ reduces equation \eqref{eq:main} to the parabolic--elliptic Keller--Segel model with a reaction term:
\begin{align*}
    \begin{cases}
        u_t = du_{xx} - k(uv_x)_x + f(u), \\
        0 = \beta^2v_{xx} + u - v.
    \end{cases}
\end{align*}
Therefore, the stationary problem coincides with that of a chemotaxis--growth model.
For analyses of chemotaxis--growth models, we refer to \cite{KOST12, TW07} and the references therein.

\begin{figure}[!tb]
\centering
\includegraphics[width=\linewidth]{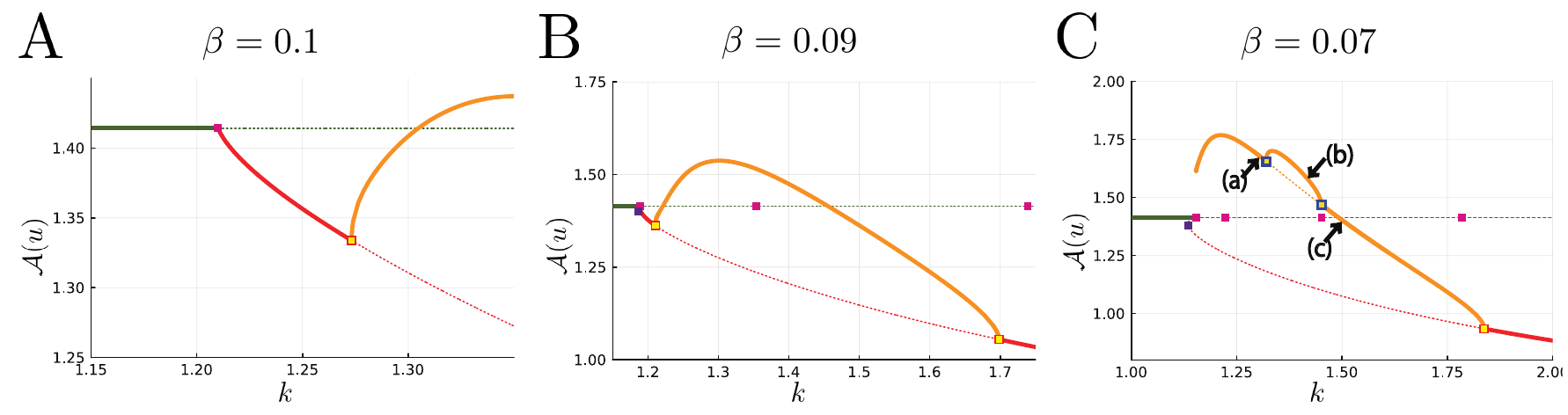}
  \caption{\small{
      Bifurcation diagrams for the exponential kernel \eqref{exp-ker}.
      See Figure~\ref{fig:gaussian} for the correspondence between branch colors and solution types.
      We set $d = c = 1.0$ and $L = 2.0$.
      A: Bifurcation diagram for $\beta = 0.1$.
      B: Bifurcation diagram for $\beta = 0.09$.
      C: Bifurcation diagram for $\beta = 0.07$.
      The bifurcation on the periodic branch is a period-doubling bifurcation.
  }}
  \label{fig:chemo}
\end{figure}

Figure~\ref{fig:chemo} shows the bifurcation diagrams for three different values of $\beta$.
As in the previous examples, we choose the parameters so that the first Fourier mode is the first mode to bifurcate.
Unlike the previous examples, even in the supercritical case, the branch of nonhomogeneous stationary solutions loses stability through a Hopf bifurcation (Figure~\ref{fig:chemo}~A).

\textbf{Multistability.}\quad
In the cases with subcritical bifurcations (Figure~\ref{fig:chemo}~B--C), the structure of the stationary solutions is qualitatively similar to those for the Gaussian and top-hat kernels.
For $\beta = 0.09$, where two Hopf bifurcations occur, there is a very narrow parameter region in which the homogeneous and nonhomogeneous stationary solutions are both stable (Figure~\ref{fig:chemo}~B).
The interval in $k$ between the primary Turing bifurcation and the saddle-node bifurcation is very narrow, and the Hopf bifurcations on this branch occur in the parameter region where the homogeneous stationary solution is unstable.
For $\beta = 0.07$, we again observe a parameter region with multistability between the homogeneous stationary solution and periodic solutions (Figure~\ref{fig:chemo}~C).

\begin{figure}[!tb]
    \centering
    \includegraphics[width=\linewidth]{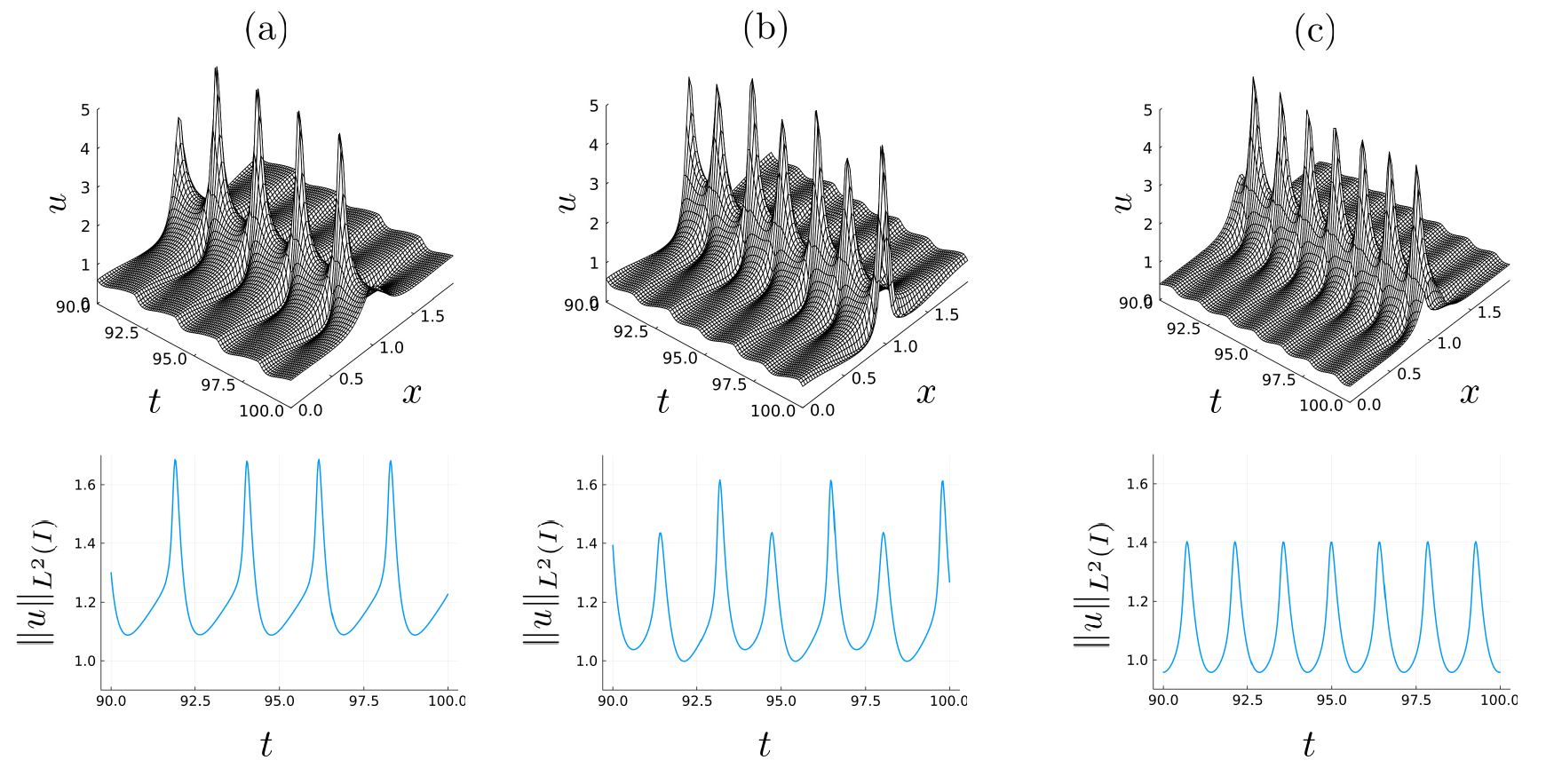}
    \caption{\small{
        Space--time plots (top row) and time evolution of the norm (bottom row) of periodic solutions for the exponential kernel \eqref{exp-ker} with $\beta = 0.07$.
        From left to right, the periodic solutions correspond to (a) $k=1.3$, (b) $k=1.4$, and (c) $k=1.5$ in Figure~\ref{fig:chemo}.
    }}
    \label{fig:chemo_po}
\end{figure}

\textbf{Period-doubling bifurcation.}\quad
For $\beta = 0.07$, where only one Hopf bifurcation occurs, a period-doubling bifurcation is observed (Figure~\ref{fig:chemo}~C).
Figure~\ref{fig:chemo_po} shows how the periodic solution changes across this bifurcation.
At points (a) and (c), corresponding to $k=1.3$ and $k=1.5$, successive oscillations have the same amplitude. 
At point (b), corresponding to $k=1.4$, large- and small-amplitude oscillations alternate.


\medskip
\subsubsection{Nonlocal-to-local approximation}
The preceding examples show how interaction localization affects bifurcation structure and multistability. 
To examine this effect from another viewpoint, we consider a formal local approximation of the convolution term. We assume that $u$ is sufficiently smooth and that $K$ has a finite second moment. 
The formal expansion is
\begin{align}\label{eq:approximate}
    (K*u)(t,x)
    \simeq
    m_0u(t,x)+\frac{m_2}{2}u_{xx}(t,x),\qquad
    m_j :=\int_{\bR}x^{j} K(x)\,dx,\ (j=0,2).
\end{align}
Since the interaction kernel is even, its first moment is zero.
This expansion was also discussed in \cite{FBC25} and used to reproduce patterns arising in the cell-sorting model of \cite{CMSTT19}.

For the rescaling considered in Theorem~\ref{thm:deformation}, decreasing $\beta$ corresponds to decreasing the second moment $m_2$ of the interaction kernel.
For the rescaled kernels considered here, the moments satisfy
\beaa
    m_0 = \int_{\bR} J(x)dx,\quad m_2 = \beta^2 \int_{\bR} x^2 J(x)dx.
\eeaa

\begin{figure}[bt]
\centering
\includegraphics[width=0.9\linewidth]{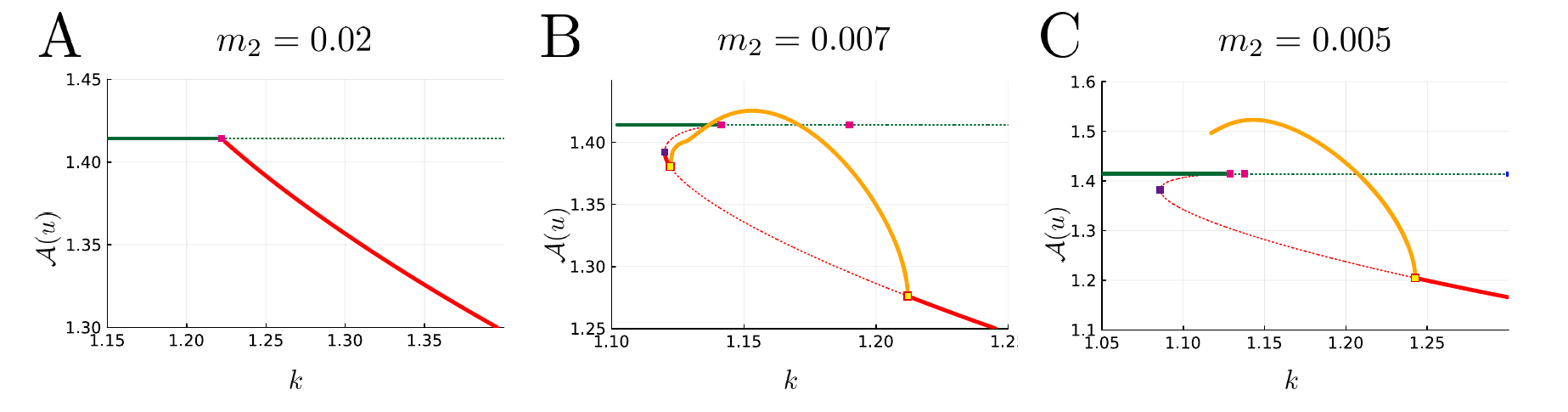}
\caption{\small{Bifurcation diagrams obtained by replacing the nonlocal term in equation \eqref{eq:main} with the approximation \eqref{eq:approximate}.
A: Bifurcation diagram for $m_2 =0.02$.
B: Bifurcation diagram for $m_2 = 0.007$.
C: Bifurcation diagram for $m_2 = 0.005$.
The other parameters are $m_0 = 1.0,\ d = c = 1.0,\ L = 2.0$.}}
\label{fig:diagram_moment}
\end{figure}

Although this approximation is outside the scope of the results in Section~\ref{sec:sta}, we performed numerical bifurcation analysis for the model obtained by replacing $K*u$ in equation~\eqref{eq:main} by the approximation \eqref{eq:approximate}, using different values of $m_2$.
The results are shown in Figure~\ref{fig:diagram_moment}.
As $m_2$ is decreased, we again observe a change from a supercritical to a subcritical bifurcation.
Time-periodic solutions are observed in both Figure~\ref{fig:diagram_moment}~B and Figure~\ref{fig:diagram_moment}~C.
In the latter case, the solution behaviour also suggests a homoclinic bifurcation.


These computations suggest that the local approximation reproduces similar bifurcation structures near the primary Turing bifurcation.
A rigorous comparison of the local and nonlocal models remains to be established.

\medskip
\subsection{Bifurcation diagrams corresponding to Case~II-A}\label{subsec:caseIIa}
In this subsection, we investigate bifurcation structure and multistability for interaction kernels corresponding to Case~II-A in Section~\ref{sec:deformation}.

\begin{figure}[!tb]
    \begin{center}
    	\includegraphics[width=\linewidth]{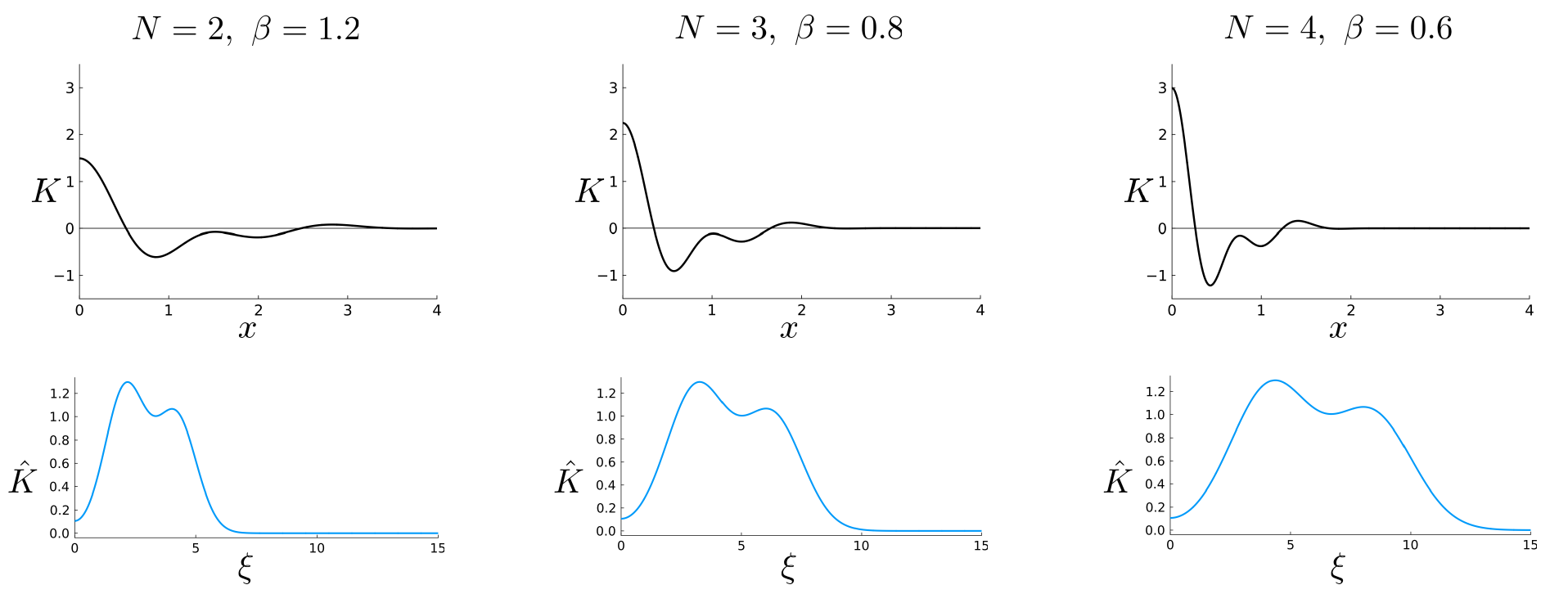}
    \end{center}
    \caption{\small{Profiles of the interaction kernel \eqref{exam:ker_subcri_case1} (top row) and its Fourier transform (bottom row).
    From left to right, the panels correspond to $(N,\beta)=(2,1.2)$, $(N,\beta)=(3,0.8)$, and $(N,\beta)=(4,0.6)$.
    The other parameters are $a_1 = 2.5,\ a_2 = 2.0,\ L = 6.0$.}}
    \label{fig:N_func}
\end{figure}

For given $N\in\bN$ and the interval length $L>0$, we use the following family of interaction kernels:
\begin{align}\label{exam:ker_subcri_case1}
    K(x) 
    = 
    a_1G(x;\beta )\cos(\xi_N x) 
    + a_2G(x;\beta )\cos(\xi_{2N}x),\quad 
    G(x;\beta) 
    = 
    \frac{1}{\sqrt{2\pi}\beta}\exp\left( -\frac{x^2}{2\beta^2}\right). 
\end{align}
The Fourier transform of this interaction kernel is given by
\beaa
    \hat{K}(\xi) 
    &=& 
    \frac{a_1\hat{G}(\xi - \xi_N;\beta)+a_1\hat{G}(\xi + \xi_N;\beta) 
    + a_2\hat{G}(\xi - \xi_{2N};\beta)+a_2\hat{G}(\xi + \xi_{2N};\beta)}{2}, \\
    \hat{G}(\xi;\beta) 
    &=& 
    \exp\left(-\frac{\beta^2 \xi^2}{2}\right).
\eeaa
We set
\beaa
    L = 6.0,\ d = c = 1.0,\ a_1 = 2.5,\ a_2 = 2.0.
\eeaa
We consider the parameter pairs $(N,\beta)=(2,1.2)$, $(3,0.8)$, and $(4,0.6)$.
For each pair, mode $N$ is the first to lose stability.
With $a_1$, $a_2$, and $L$ fixed, $R$ depends only on $\beta N$.
Since $\beta N=2.4$ for all three pairs, they have the same value $R\simeq0.8147940168$. 
Figure~\ref{fig:N_func} shows the kernels and their Fourier transforms. 
These kernels produce both attractive and repulsive interactions, depending on distance.

\begin{figure}[!tb]
  \centering
    \includegraphics[width=\linewidth]{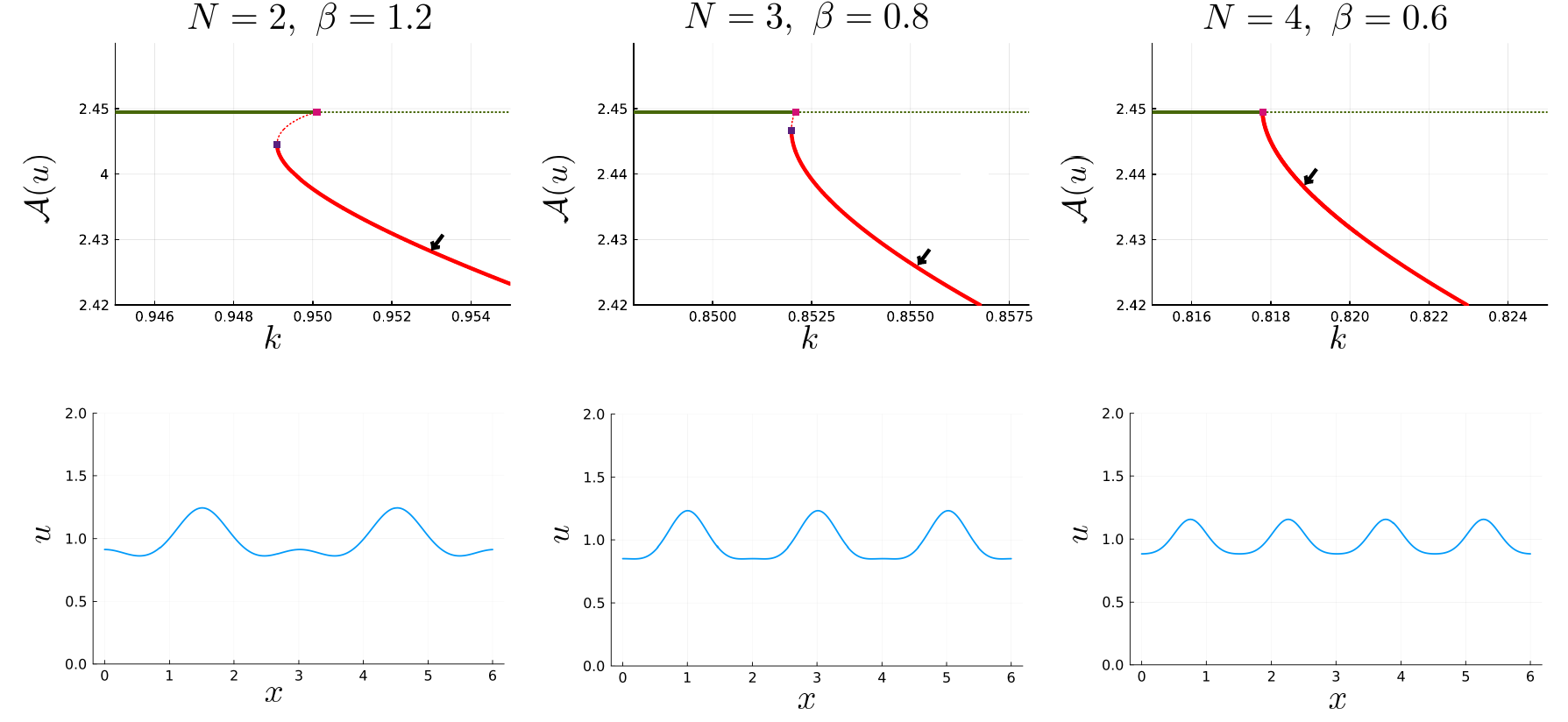}
  \caption{\small{Bifurcation diagrams obtained using the interaction kernels in Figure~\ref{fig:N_func} (top row) and spatial profiles of stable nonhomogeneous stationary solutions (bottom row).
  The branch colors and corresponding solution types are the same as in Figure~\ref{fig:gaussian}.
  From left to right, the panels correspond to $(N,\beta)=(2,1.2)$, $(N,\beta)=(3,0.8)$, and $(N,\beta)=(4,0.6)$.
  The stable nonhomogeneous stationary solutions are shown at $k=0.953031$, $0.8551655$, and $0.8189258$, respectively, as indicated by the arrows in the corresponding bifurcation diagrams.
  The other parameters are $a_1 = 2.5,\ a_2 = 2.0,\ L = 6.0$.}}
  \label{fig:diagram_case1}
\end{figure}

For kernels in Case II-A, subcriticality requires $\alpha=d\xi_N^2$ to lie in the bounded interval where $F(\alpha)>0$.
Table 1 shows that inequality \eqref{F(x)>0} holds for $N=2$ and $N=3$, but not for $N=4$.
Accordingly, the bifurcation is subcritical for $N=2,3$ and supercritical for $N=4$.
Figure~\ref{fig:diagram_case1} shows the corresponding bifurcation diagrams.


In the two subcritical cases, the stable homogeneous stationary solution coexists with stable nonhomogeneous stationary solutions over a range of $k$. 
No time-periodic solutions were detected along the computed branches near the primary Turing bifurcation.

\begin{table}[b]
 \begin{center}
   \caption{Verification of inequality \eqref{F(x)>0}}\label{table:caseII-a}
  \begin{tabular}{|c|c|c|}\hline
    Parameter  & $\alpha = d\xi^2_{N}$ & \eqref{F(x)>0}  \\ \hline
    $N=2$, $\beta=1.2$  & 4.38649084493 &true \\ 
    $N=3$, $\beta=0.8$  & 9.86960440109 &true \\
    $N=4$, $\beta=0.6$  & 17.5459633797 &false \\ \hline
  \end{tabular}
 \end{center}
\end{table}

\medskip

\subsection{Bifurcation diagrams corresponding to Case~II-B}\label{subsec:caseIIb}
\begin{figure}[!tb]
  \centering
    \includegraphics[width=\linewidth]{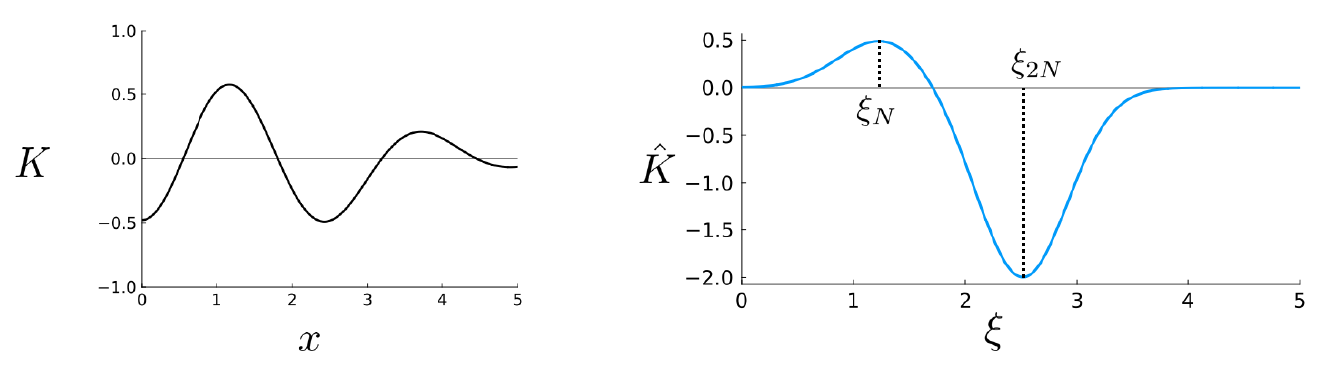}
  \caption{\small{Interaction kernel \eqref{exam:ker_subcri_case1} with $a_1 = 1.0,\ a_2 = -4.0,\ \beta = 2.5$ (left) and its Fourier transform (right).}}
  \label{fig:subcri_plot}
\end{figure}

Finally, we consider a setting corresponding to Case~II-B in Section~\ref{sec:deformation}.
We again use the kernel \eqref{exam:ker_subcri_case1} with
\beaa
    N=2,\ L=10.0,\ a_1=1.0,\ a_2=-4.0,\ \beta=2.5.
\eeaa
For these parameter values, $\hat{K}(\xi_{N})>0$, whereas $\hat{K}(\xi_{2N})<0$.
Unlike the interaction kernels considered in the previous section, this interaction kernel has a repulsive effect at short distances (Figure~\ref{fig:subcri_plot}).


The ratio
\beaa
    R \simeq -4.11 < -\dfrac{1}{6}
\eeaa
places this kernel in Case~II-B.
We set $d=0.05$ and $c=1.0$.
For these values, $\alpha=d\xi_N^2\simeq0.07896$ satisfies $F(\alpha)>0$, so the Turing bifurcation is subcritical.

\begin{figure}[bt]
\begin{center}
	\includegraphics[width=\linewidth]{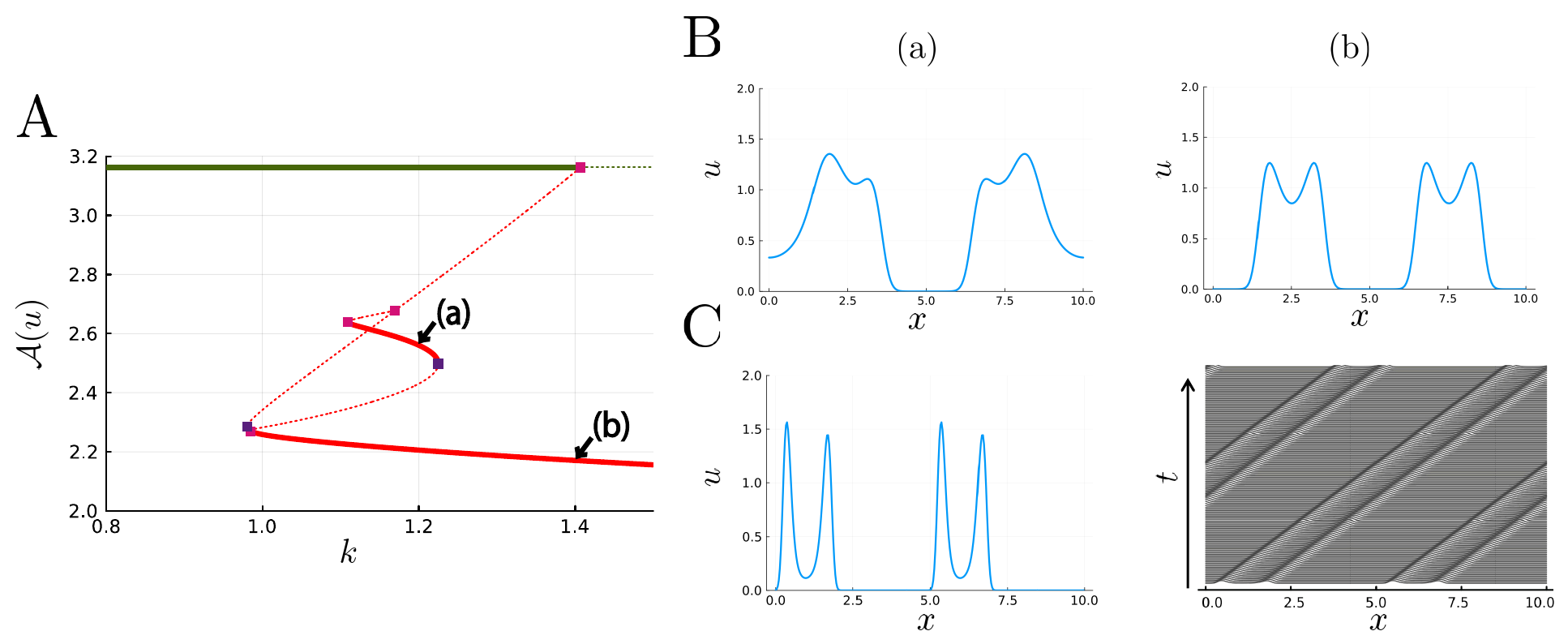}
\end{center}
\caption{\small{
The branch colors and corresponding solution types are the same as in Figure~\ref{fig:gaussian}.
A: Bifurcation diagram for the interaction kernel \eqref{exam:ker_subcri_case1} with $a_1 = 1.0, a_2 = -4.0, \beta = 2.5$.
B: Stable stationary solutions corresponding to (a) and (b) in the bifurcation diagram.
C: Profile and surface plot of the traveling wave observed at $k = 19.0$.
The other parameters are $d = 0.05,\ c = 1.0,\ L = 10.0$.
The first unstable mode is $N=2$.
}}
\label{fig:case2}
\end{figure}

Figure~\ref{fig:case2}~A  shows a parameter interval in which the stable homogeneous stationary solution coexists with two distinct stable nonhomogeneous stationary solutions. 
This interval of multistability is wider than in the preceding numerical examples.
The branch labeled (b) is not stabilized by the saddle-node bifurcation near $k=1.0$.
Instead, it becomes stable through a nearby pitchfork bifurcation.


Figure~\ref{fig:case2}~B shows the stable nonhomogeneous stationary solutions corresponding to (a) and (b) in Figure~\ref{fig:case2}~A.
For the small diffusion coefficient used here, the profiles in Figure~\ref{fig:case2}~B have distinct high- and low-density regions.
The stationary branch indicated by (b) later loses stability near $k=19.0$.
Numerical simulations at $k=19.0$ reveal a traveling wave with an asymmetric spatial profile (Figure~\ref{fig:case2}~C).
This wave is observed far from the primary Turing bifurcation.


\bigskip
\section{Concluding remarks}\label{sec:dis}
In this study, we analyzed Turing bifurcations in a nonlocal advection--reaction--diffusion equation with logistic growth. 
We first established the existence of nonhomogeneous stationary solutions bifurcating from the positive homogeneous stationary solution. 
Center manifold reduction yielded an amplitude equation near the bifurcation point. 
Its cubic coefficient provided criteria for supercritical and subcritical bifurcations in terms of the Fourier transform of the kernel. 
We also derived a sufficient condition for subcriticality based on kernel rescaling. 
Under the assumptions of Theorem~\ref{thm:deformation}, sufficiently localized attractive interactions can produce subcritical Turing bifurcations.


Numerical bifurcation analysis revealed several forms of multistability. 
The stable homogeneous stationary solution coexisted with stable nonhomogeneous stationary solutions or stable time-periodic solutions. 
In particular, coexistence with time-periodic solutions was observed for strongly localized attractive interactions.
These results suggest that parameter changes or finite-amplitude perturbations can induce abrupt transitions from the homogeneous state to large-amplitude periodic behavior.


Multistability was also observed for selected kernels in the oscillatory family \eqref{exam:ker_subcri_case1}. 
In one example, the stable homogeneous stationary solution coexisted with two distinct stable nonhomogeneous stationary solutions. 
Traveling waves were also observed numerically far from the primary Turing bifurcation. 
However, these oscillatory kernels do not have an obvious biological or physical interpretation. 
The Gaussian and top-hat kernels illustrate how rescaling can produce subcriticality in more interpretable interaction models.
Extending this analysis to a broader range of physically interpretable kernel families is an important direction for future work. 
Such extensions may clarify the relation between interaction properties and bifurcation structure.


Another important direction is to extend the analysis to two spatial dimensions. 
Previous studies have shown that pattern types in nonlocal advection models can depend on spatial dimension \cite{JKMG23}.
Such an extension is also relevant to animal aggregation and the motion of cell populations. 
A natural question is how the relation between kernel shape and bifurcation criticality obtained here extends to this setting.


Models incorporating nonlinear diffusion and volume-filling effects have also been proposed to describe cell sorting \cite{CMSTT19}.
Although these models can reproduce experimental observations in numerical simulations, their degeneracy may lead to compactly supported solutions.
This makes standard numerical bifurcation methods difficult to apply directly. 
Developing a numerical bifurcation framework for these models could deepen our understanding of their bifurcation structures and mathematical properties. 
It could also help clarify their relation to experimentally observed phenomena.

\bigskip
\section*{Declaration of generative AI and AI-assisted technologies in the writing process}
During the preparation of this work, the authors used DeepL and ChatGPT in order to improve their English writing.
After using these tools, the authors reviewed and edited the content as needed and take full responsibility for the content of the publication.

\bigskip
\section*{Acknowledgments}
MN is supported by JSPS KAKENHI Grant Number JP25K00918 and MEXT KAKENHI Grant Number JP23H04936.
HI is supported by JSPS KAKENHI Grant Numbers JP23K13013, JP24H00188 and JP25K07142.

\bigskip
\section*{Declaration of Competing Interest}
The authors declare that they have no known competing financial interests or personal relationships that could have appeared to influence the work reported in this paper.

\appendix

\bigskip
\section{Proof of Proposition~\ref{prop:cr}}\label{sec:A}
\renewcommand{\theequation}{A.\arabic{equation}}
In this section, we prove Proposition~\ref{prop:cr} by applying the Crandall--Rabinowitz bifurcation theorem.
We take $k>0$ as the bifurcation parameter and use the function spaces $\mathcal{X}$ and $\mathcal{Y}$ defined in \eqref{eq:function-spaces}.
Writing $u=1+w$, we define the operators from $\mathcal{X}$ to $\mathcal{Y}$ by
$$
    \mathcal{L}_1w := dw_{xx} -cw,\quad
    \mathcal{L}_2 w := K*w_{xx},\quad
    \mathcal{N}_1(w) := -cw^2, \quad
    \mathcal{N}_2(w) := (w(K*w)_x)_x.
$$
As stated in the main text, $\mathcal{L}_1,  \mathcal{L}_2: \mathcal{X} \to \mathcal{Y}$, and the nonlinear terms are well defined.
Then the stationary problem for equation \eqref{eq:main}, expanded around $u = 1$, can be written as
$$
    \mathcal{H}(k, w; d, c) 
    := \mathcal{L}_1w + \mathcal{N}_1(w) - k(\mathcal{L}_2 w + \mathcal{N}_2(w)) = 0,
$$
where $\mathcal{H} : \bR \times \mathcal{X} \to \mathcal{Y}$.
For every $k$,
$$
\mathcal{H}(k,0;d,c)=0.
$$
Moreover, the nonlinear terms are of class $C^2$ and, in particular,
$\mathcal{N}_i(0) = 0$ and $D\mathcal{N}_i(0) = 0$ for $i=1,2$.
Here, $D\mathcal{N}_i(w)$ denotes the Fr\'{e}chet derivative of $\mathcal{N}_i$ at $w$.
We first establish the following properties of the linear operators.

\begin{lemma}\label{lem:linearop}
    Suppose that Assumption~\ref{assumption:lam} holds.
    Then the linear operators $\mathcal{L}_1,\mathcal{L}_2:\mathcal{X}\to\mathcal{Y}$ have the following properties:
    \begin{enumerate}
        \item[(i)] The kernel of $\mathcal{L}_1-k^*\mathcal{L}_2$ is given by
        \[
            \operatorname{Ker}(\mathcal{L}_1-k^*\mathcal{L}_2)
            =\operatorname{span}\{\,\cos(\xi_N x)\,\};
        \]
        \item[(ii)] The range of $\mathcal{L}_1-k^*\mathcal{L}_2$ is given by
        \[
            \operatorname{Range}(\mathcal{L}_1-k^*\mathcal{L}_2)
            =\left\{
                v\in\mathcal{Y}
                \mid (v,\cos(\xi_N x))_{L^2(I)}=0
            \right\};
        \]
        \item[(iii)] The following transversality condition holds:
        \[
            \int_I \bigl(\mathcal{L}_2\cos(\xi_N x)\bigr) \cos(\xi_N x)\,dx\neq0.
        \]
    \end{enumerate}
\end{lemma}
\begin{proof}[Proof of Lemma~\ref{lem:linearop}]
    \textit{(i)} Every $w\in\mathcal{X}$ admits a cosine expansion
    \be\label{cos}
        w(x)=\sum_{n\ge0}a_n\cos(\xi_n x),
    \ee
    and
    \beaa
        (\mathcal{L}_1-k^*\mathcal{L}_2)w(x)
        =\sum_{n\ge0}\lambda_n^*a_n\cos(\xi_n x).
    \eeaa
    Since $\lambda_0^*=-c<0$, Assumption~\ref{assumption:lam} implies that
    $\lambda_n^*=0$ if and only if $n=N$.
    Hence,
    \beaa
        \operatorname{Ker}(\mathcal{L}_1-k^*\mathcal{L}_2) =\operatorname{span}\{\,\cos(\xi_N x)\,\}.
    \eeaa

    \medskip
    \noindent\textit{(ii)}
    The expansion in (i) and the identity $\lambda_N^*=0$ give
    \beaa
        \operatorname{Range}(\mathcal{L}_1-k^*\mathcal{L}_2) \subset
        \left\{  v\in\mathcal{Y} \mid (v,\cos(\xi_N x))_{L^2(I)}=0 \right\}.
    \eeaa
    To prove the reverse inclusion, let $v\in\mathcal{Y}$ satisfy
    $(v,\cos(\xi_N x))_{L^2(I)}=0$.
    Write
    \beaa
        v(x)=\sum_{n\ge0}b_n\cos(\xi_n x).
    \eeaa
    By orthogonality, $b_N=0$.
    Define
    \beaa
        a_N=0,\qquad
        a_n=\frac{b_n}{\lambda_n^*}
        \quad(n\ge0,\ n\neq N)
    \eeaa
    and define $w$ by the cosine expansion \eqref{cos}.
    We show that $w\in\mathcal{X}$.

    By the Riemann--Lebesgue lemma, $\hat{K}(\xi)\to 0$ as $|\xi|\to\infty$.
    Since
    \beaa
        \lambda_n^* =\bigl(-d+k^*\hat{K}(\xi_n)\bigr)\xi_n^2-c,
    \eeaa
    it follows that
    \beaa
        |\lambda_n^*|\ge\frac{d}{2}\xi_n^2
    \eeaa
    for all sufficiently large $n$.
    The finitely many remaining eigenvalues with $n\neq N$ are nonzero.
    Therefore,
    \beaa
        C_0:=
        \sup_{\substack{n\ge0\\n\neq N}}
        \frac{1+\xi_n^4}{|\lambda_n^*|^2}
        <\infty.
    \eeaa
    Since $v\in L^2(I)$, we obtain
    \begin{align*}
        \sum_{n\ge0}(1+\xi_n^4)|a_n|^2
        =
        \sum_{\substack{n\ge0\\n\neq N}}
        \frac{1+\xi_n^4}{|\lambda_n^*|^2}|b_n|^2
        \le C_0\sum_{n\ge0}|b_n|^2
        <\infty.
    \end{align*}
    Thus, the cosine expansion \eqref{cos} converges in $H^2_{\mathrm{per}}(I)$ and defines an element of $\mathcal{X}$.
    By construction,
    \beaa
        (\mathcal{L}_1-k^*\mathcal{L}_2)w=v.
    \eeaa
    This proves the reverse inclusion and hence the claimed
    characterization of the range.

    Define the linear functional $y^*:\mathcal{Y}\to\mathbb{R}$ by
    \beaa
        y^*(v):=(v,\cos(\xi_N x))_{L^2(I)}.
    \eeaa
    By the Cauchy--Schwarz inequality, $y^*$ is continuous.
    It is also nonzero, since
    \beaa
        y^*(\cos(\xi_N x)) =\int_I\cos^2(\xi_N x)\,dx =\frac{L}{2}>0.
    \eeaa
    Consequently,
    \beaa
        \operatorname{Range}(\mathcal{L}_1-k^*\mathcal{L}_2)
        =\operatorname{Ker}(y^*)
    \eeaa
    is a closed subspace of $\mathcal{Y}$ with
    \beaa
        \operatorname{codim}
        \operatorname{Range}(\mathcal{L}_1-k^*\mathcal{L}_2)=1.
    \eeaa

    \medskip
    \noindent\textit{(iii)}
    Since $K$ is even,
    $K*\cos(\xi_N x)=\hat{K}(\xi_N)\cos(\xi_N x)$ holds.
    Using $\mathcal{L}_2w=K*w_{xx}$, we obtain
    \begin{align*}
        \int_I
        \bigl(\mathcal{L}_2\cos(\xi_N x)\bigr)
        \cos(\xi_N x)\,dx
        =-\frac{L\xi_N^2\hat{K}(\xi_N)}{2}\neq0,
    \end{align*}
    where the nonvanising follows from $\lambda_N^*=0$.
    This proves (iii) and completes the proof.
\end{proof}

We can now apply the Crandall--Rabinowitz bifurcation theorem \cite{CR71} to
$\mathcal H(k, w;d, c) = 0$ at $(k^*, 0)$.
This proves Proposition~\ref{prop:cr}.

\bigskip

\section{Local and global existence}\label{appendix:exi}
\renewcommand{\theequation}{B.\arabic{equation}}
In this appendix, we prove Propositions~\ref{prop:local_ex} and \ref{prop:global_ex}.
We first construct a local mild solution by the contraction mapping principle and then show that it is a classical solution by Schauder estimates.
Next, on any time interval $[0,T]$ on which the solution exists, we derive $L^p$ and $L^\infty$ estimates that are independent of $T$.
Finally, using these uniform estimates together with the local existence result, we extend the solution to any finite time.

\subsection{Proof of Proposition~\ref{prop:local_ex}}
We first prove the local existence result stated in Proposition~\ref{prop:local_ex}.
Define the function space
$$
    \mathcal{X}_T
    :=
    \left\{\, u \in C([0,T]; C_{\mathrm{per}}(I)) \ \middle|\ 
        u(t,\cdot)\in C^1_{\mathrm{per}}(I) \ (t>0) ,\ 
        \sup_{0 < t \le T}\sqrt{t}\|u_x(t)\|_{L^{\infty}(I)}
        <
        \infty\,
    \right\}
$$
with the norm
$$
    \| u \|_{\mathcal{X}_T}
    :=
    \max_{0\le t \le T}\|u(t)\|_{L^{\infty}(I)}
    +
    \sup_{0 < t \le T}\sqrt{t}\|u_x(t)\|_{L^{\infty}(I)}.
$$
With this norm, $\mathcal{X}_T$ is a Banach space.
We also define the periodic heat kernel by
\begin{align*}
    H_d^L(t,x)
    :=
    \sum_{j\in\mathbb{Z}}H_d(t,x-jL),
    \qquad
    H_d(t,x)
    :=
    \frac{1}{2\sqrt{\pi dt}}e^{-x^2/(4dt)}.
\end{align*}


\begin{definition}
    Let $\phi \in C_{\mathrm{per}}(I)$ be an initial value.
    A function $u \in \mathcal{X}_T$ is called a mild solution of equation~\eqref{eq:main} on $[0,T]$ if it satisfies
    \begin{align}
        u(t)
        ={}&
        H_d^L(t)*\phi
        -
        k\int_0^t\,\partial_xH_d^L(t-s)*\bigl[u(K'*u)\bigr](s)\,ds 
        +
        \int_0^t\,H_d^L(t-s)*f(u(s))\,ds.
        \label{eq:mild_appendix}
    \end{align}
\end{definition}


Since $f$ is locally Lipschitz continuous, for every $R>0$, there exists $L_R>0$ such that
$$
    |f(u)-f(v)|\le L_R|u-v|
$$
whenever $|u|,|v|\le R$.
We also use $f(0)=0$ below.
Set
$$
    C:=\frac{1}{\sqrt{\pi}}.
$$


\begin{lemma}\label{lem:contraction}
    Suppose that $K\in W^{1,1}(\mathbb{R})$.
    Let 
    \begin{align}
        R
        >
        2\left(\frac{C}{\sqrt d}+1\right)\|\phi\|_{L^\infty(I)},
        \label{eq:R_condition}
    \end{align}
    and $T$ be a constant satisfying
    \begin{align}
        \frac{2kC(1+\pi)}{\sqrt d}R\|K'\|_{L^1(\bR)}\sqrt T
        +
        \left(1+\frac{2C}{\sqrt d}\right)L_RT
        <
        \frac{1}{2}.\label{eq:T_condition}
    \end{align}
    Define
    $$
        B_R(\mathcal X_T)
        :=
        \{\,u\in\mathcal X_T\mid \|u\|_{\mathcal X_T}\le R\,\}.
    $$
    Then, equation~\eqref{eq:main} admits a unique mild solution
    $u\in B_R(\mathcal X_T)$ on $[0,T]$.
\end{lemma}

\begin{proof}
    Let $\Phi$ be the map defined by the right-hand side of \eqref{eq:mild_appendix}.
    We show that $\Phi$ maps $B_R(\mathcal X_T)$ into itself and is a contraction.
    Set $F(u):=u(K'*u)$.
    For $u\in B_R(\mathcal X_T)$, Young's inequality gives
    \begin{align*}
        \|F(u(s))\|_{L^\infty(I)}
        \le
        R^2\|K'\|_{L^1(\mathbb R)},\qquad
        \|\partial_xF(u(s))\|_{L^\infty(I)}
        \le
        2R^2\|K'\|_{L^1(\mathbb R)}s^{-1/2}.
    \end{align*}
    Here, we have used
    $$
        \partial_xF(u)
        =
        u_x(K'*u)+u(K'*u_x).
    $$
    In particular, the identity $\partial_x(K'*u)=K'*u_x$ follows from the periodicity of $u$ and the assumption $K'\in L^1(\mathbb R)$.
    The periodic heat kernel satisfies
    \begin{align*}
        \|H_d^L(t)\|_{L^1(I)}=1,\qquad
        \|\partial_xH_d^L(t)\|_{L^1(I)}
        \le
        \frac{C}{\sqrt{dt}}.
    \end{align*}
    Standard properties of the heat semigroup also show that $\Phi[u]\in\mathcal X_T$ for every $u\in B_R(\mathcal X_T)$.
    Therefore,
    \begin{align*}
        \|\Phi[u]\|_{\mathcal X_T}
        \le{}
        \left(\frac{C}{\sqrt d}+1\right)\|\phi\|_{L^\infty(I)} 
        +
        \frac{2kC(1+\pi)}{\sqrt d}R^2\|K'\|_{L^1(\mathbb R)}\sqrt T
        +
        \left(1+\frac{2C}{\sqrt d}\right)L_RRT.
    \end{align*}
    By \eqref{eq:R_condition} and \eqref{eq:T_condition},
    \begin{gather*}
        \left(\frac{C}{\sqrt d}+1\right)\|\phi\|_{L^\infty(I)}
        <
        \frac{R}{2},\quad
        \frac{2kC(1+\pi)}{\sqrt d}R^2\|K'\|_{L^1(\mathbb R)}\sqrt T
        +
        \left(1+\frac{2C}{\sqrt d}\right)L_RRT
        <
        \frac{R}{2}.
    \end{gather*}
    Hence,
    $\|\Phi[u]\|_{\mathcal X_T}<R$, and therefore
    $\Phi(B_R(\mathcal X_T))\subset B_R(\mathcal X_T)$.

    Next, let $u,v\in B_R(\mathcal X_T)$.
    Similarly,
    \begin{align*}
        \|F(u(s))-F(v(s))\|_{L^\infty(I)}
        &\le
        2R\|K'\|_{L^1(\mathbb R)}\|u-v\|_{\mathcal X_T},\\
        \|\partial_xF(u(s))-\partial_xF(v(s))\|_{L^\infty(I)}
        &\le
        4R\|K'\|_{L^1(\mathbb R)}s^{-1/2}\|u-v\|_{\mathcal X_T}.
    \end{align*}
    It follows that
    \begin{align}
        \|\Phi[u]-\Phi[v]\|_{\mathcal X_T}
        \le{}&
        \left\{
        \frac{4kC(1+\pi)}{\sqrt d}R\|K'\|_{L^1(\mathbb R)}\sqrt T
        +
        \left(1+\frac{2C}{\sqrt d}\right)L_RT
        \right\}
        \|u-v\|_{\mathcal X_T}.
        \label{eq:Phi_contraction_bound}
    \end{align}
    By the assumption on $T$,
    \begin{align*}
        \frac{4kC(1+\pi)}{\sqrt d}R\|K'\|_{L^1(\mathbb R)}\sqrt T
        +
        \left(1+\frac{2C}{\sqrt d}\right)L_RT<1.
    \end{align*}
    Hence, $\Phi$ is a contraction on $B_R(\mathcal X_T)$.
    The contraction mapping principle therefore yields a unique fixed point $u\in B_R(\mathcal X_T)$, which is a mild solution of equation~\eqref{eq:main}.
\end{proof}




We next show that the mild solution obtained above is a classical solution.

\begin{lemma}\label{lem:schauder}
    Let $u\in\mathcal X_T$ be the mild solution obtained above.
    Then, for every $\tau\in(0,T)$ and $\alpha\in(0,1)$,
    $$
        u\in C^{1+\alpha/2,\,2+\alpha}([\tau,T]\times I).
    $$
\end{lemma}

\begin{proof}
    By standard parabolic smoothing applied to the mild formulation, for every $\tau>0$,
    $$
        u\in C^{\alpha/2,\,1+\alpha}([\tau,T]\times I).
    $$
    Since $K\in W^{1,1}(\mathbb R)$, we have $K'\in L^1(\mathbb R)$, and convolution with $K'$ is bounded on $C^\alpha_{\mathrm{per}}(I)$.
    Hence,
    \begin{align*}
        [u(K'*u)]_x
        =
        u_x(K'*u)+u(K'*u_x)
    \end{align*}
    belongs to $C^{\alpha/2,\,\alpha}([\tau,T]\times I)$.
    The same holds for $f(u)$.
    Therefore, regarding equation~\eqref{eq:main} as a nonhomogeneous linear parabolic equation, the conclusion follows from the Schauder estimates \cite[Theorem~9.1.2]{Krylov96}.
\end{proof}

Finally, we prove uniqueness of classical solutions.
Let $u_1$ and $u_2$ be two classical solutions on $[0,T]$ in the class stated in Proposition~\ref{prop:local_ex}, with $u_j(0)=\phi$ for $j=1,2$.
Both satisfy the mild formulation \eqref{eq:mild_appendix}.
Set
\beaa
     M := 1+\sup_{0\le t\le T} \max\left\{ \|u_1(t)\|_{L^\infty(I)}, \|u_2(t)\|_{L^\infty(I)} \right\}.
\eeaa
For $0<h\le T$, define
\beaa
    D_h := \sup_{0\le s\le h} \|u_1(s)-u_2(s)\|_{L^\infty(I)}.
\eeaa
Using \eqref{eq:mild_appendix}, the heat kernel estimates, and the local
Lipschitz continuity of $f$, we obtain
\beaa
    D_h \le \left( \frac{4kCM\|K'\|_{L^1(\mathbb R)}}{\sqrt d}\sqrt h +L_Mh \right)D_h.
\eeaa
Choose $h>0$ sufficiently small that the factor in parentheses is less than one. Then $D_h=0$, so $u_1=u_2$ on $[0,h]$.
Repeating the same argument on successive time intervals of length at most $h$, with the same bound $M$, proves that $u_1=u_2$ on $[0,T]$.
This completes the proof of Proposition~\ref{prop:local_ex}.

\medskip
\subsection{Proof of Proposition~\ref{prop:global_ex}}
We next derive the estimates needed to prove global existence.
In what follows, we assume the hypotheses of Proposition~\ref{prop:global_ex}, in particular,
$K\in W^{2,1}(\bR)$, $f(u)=cu(1-u)$, $c>0$, and $\phi\ge0$.
Let $u$ be a classical solution on a time interval $[0,T]$.


We first establish the nonnegativity of the solution.

\begin{lemma}\label{lem:positive}
    If $\phi\ge0$, then $u(t,x)\ge 0$ for all $(t,x)\in[0,T]\times I$.
    If, in addition, $\phi\not\equiv 0$, then $u(t,x)> 0$ for all $(t,x)\in(0,T]\times I$.
\end{lemma}
\begin{proof}
    Since $K\in W^{2,1}(\bR)$, equation~\eqref{eq:main} can be written as
    $$
        u_t-du_{xx}+k(K'*u)u_x
        =
        u\bigl[-k(K''*u)+c(1-u)\bigr].
    $$
    The coefficients are bounded on $[0,T]\times I$. 
    The maximum principle (e.g. \cite{PW84}) on the periodic domain therefore yields $u\ge0$.
    If $\phi\not\equiv0$, the strong maximum principle implies that $u(t,x)>0$ for every $t\in(0,T]$ and $x\in I$.
\end{proof}

In the next lemma, we first derive an $L^1$ estimate and then use it to obtain $L^p$ estimates for every finite $p$.
The proof follows the argument of \cite[Lemma~2.2]{TW07}.

\begin{lemma}\label{lem:Lp_bound}
    For every $p\in[1,\infty)$, there exists
    $C=C(\|\phi\|_{L^\infty(I)},p)>0$ such that
    $$
        \|u(t)\|_{L^p(I)}\le C
    $$
    for all $t\in[0,T]$.
    Here, $C$ is independent of $T$.
\end{lemma}

\begin{proof}
    We first consider the case $p=1$.
    By Lemma~\ref{lem:positive}, we have $u\ge0$, and we set
    $$
        M(t):=\int_Iu(t)\,dx=\|u(t)\|_{L^1(I)}.
    $$
    Using the periodic boundary conditions and H\"{o}lder's inequality, we obtain
    \begin{align*}
        M'(t) = \int_I cu(1-u)\,dx \le  cM(t)-\frac{c}{|I|}M(t)^2.
    \end{align*}
    Hence,
    \begin{align}
        M(t)
        \le
        M_0
        :=
        \max\{\,\|\phi\|_{L^1(I)},\ |I|\,\}.
        \label{eq:L1_bound}
    \end{align}


    Next, let $p\ge2$.
    Multiplying equation~\eqref{eq:main} by $u^{p-1}$ and integrating over $I$, we obtain
    \begin{align*}
        \frac{1}{p}\frac{d}{dt}\|u\|_{L^p(I)}^p
        =& -\frac{4d(p-1)}{p^2}\|\partial_xu^{p/2}\|_{L^2(I)}^2+
        \frac{2k(p-1)}{p}
        \int_I(\partial_xu^{p/2})u^{p/2}(K'*u)\,dx\\
        &\qquad\qquad+
        c\|u\|_{L^p(I)}^p
        -
        c\|u\|_{L^{p+1}(I)}^{p+1}.
    \end{align*}
    Applying Young's inequality to the nonlocal term gives
    \begin{align*}
        \frac{2k(p-1)}{p}
        \int_I(\partial_xu^{p/2})u^{p/2}(K'*u)\,dx
        \le
        \frac{d(p-1)}{p^2}\|\partial_xu^{p/2}\|_{L^2(I)}^2
        +
        \frac{k^2(p-1)}{d}
        \|u^{p/2}(K'*u)\|_{L^2(I)}^2.
    \end{align*}


    Let
    $$
        K'_L(x):=\sum_{j\in\mathbb Z}K'(x+jL)
    $$
    denote the periodization of $K'$.
    Then the periodic convolution $K'*u$ can be written using $K'_L$, and
    \begin{align*}
        \|u^{p/2}(K'*u)\|_{L^2(I)}^2
        \le
        \|u\|_{L^p(I)}^p\|K'*u\|_{L^\infty(I)}^2
        \le
        \|K'_L\|_{L^\infty(I)}^2M_0^2\|u\|_{L^p(I)}^p.
    \end{align*}
    Since $K\in W^{2,1}(\mathbb R)$, we have
    \begin{align*}
        \|K'_L\|_{L^\infty(I)}
        \le
        \frac{1}{L}\|K'_L\|_{L^1(I)}
        +
        \|(K'_L)'\|_{L^1(I)}
        \le
        \frac{1}{L}\|K'\|_{L^1(\mathbb R)}
        +
        \|K''\|_{L^1(\mathbb R)}
        <
        \infty.
    \end{align*}
    Combining these estimates, we obtain
    \begin{align}
        \frac{1}{p}\frac{d}{dt}\|u\|_{L^p(I)}^p
        &\le
        \left(
            c+\frac{k^2(p-1)}{d}\|K'_L\|_{L^\infty(I)}^2M_0^2
        \right)
        \|u\|_{L^p(I)}^p
        -
        c\|u\|_{L^{p+1}(I)}^{p+1}.
        \label{ineq:Lp}
    \end{align}
    Set
    $$
        P
        :=
        c+\frac{k^2(p-1)}{d}\|K'_L\|_{L^\infty(I)}^2M_0^2.
    $$
    Taking $\alpha=(p+1)/p$ and $\beta=p+1$, so that $1/\alpha+1/\beta=1$, Young's inequality gives
    \begin{align*}
        Pu^p
        \le
        \frac{c}{2}u^{p+1}
        +
        \frac{1}{\beta}
        \left(\frac{\alpha c}{2}\right)^{-p}P^{p+1}.
    \end{align*}
    Denote the second term on the right-hand side by $C_p$.
    Then \eqref{ineq:Lp} yields
    \begin{align*}
        \frac{1}{p}\frac{d}{dt}\|u\|_{L^p(I)}^p
        \le
        -\frac{c}{2}\|u\|_{L^{p+1}(I)}^{p+1}
        +
        C_p|I|.
    \end{align*}
    Furthermore, H\"older's inequality gives
    \begin{align*}
        \|u\|_{L^{p+1}(I)}^{p+1}
        \ge
        |I|^{-1/p}\|u\|_{L^p(I)}^{p+1}.
    \end{align*}
    Thus, setting $Y(t):=\|u(t)\|_{L^p(I)}^p$, we obtain
    \begin{align*}
        Y'(t)
        \le
        -a_pY(t)^{(p+1)/p}+b_p,
    \end{align*}
    where
    $$
        a_p:=\frac{pc}{2|I|^{1/p}},
        \qquad
        b_p:=pC_p|I|.
    $$
    Therefore,
    \begin{align*}
        Y(t)
        \le
        \max\left\{
            Y(0),
            \left(\frac{b_p}{a_p}\right)^{p/(p+1)}
        \right\}.
    \end{align*}
    Hence,
    \begin{align}
        \|u(t)\|_{L^p(I)}
        \le
        \max\left\{
            \|\phi\|_{L^p(I)},
            \left(\frac{b_p}{a_p}\right)^{1/(p+1)}
        \right\}.
        \label{eq:Lp_bound_pge2}
    \end{align}


    Finally, consider the case $1<p<2$.
    Let $\theta\in(0,1)$ be defined by
    $$
        \frac{1}{p}
        =
        \theta+\frac{1-\theta}{2}.
    $$
    By H\"older's inequality,
    \begin{align*}
        \|u(t)\|_{L^p(I)}
        \le
        \|u(t)\|_{L^1(I)}^\theta
        \|u(t)\|_{L^2(I)}^{1-\theta}.
    \end{align*}
    The conclusion then follows from \eqref{eq:L1_bound} and the case $p=2$ of \eqref{eq:Lp_bound_pge2}.
\end{proof}


We next use the estimates obtained in Lemma~\ref{lem:Lp_bound} to derive an $L^\infty$ estimate.
The periodic heat semigroup satisfies the following estimates, analogous to those used in \cite[Lemma~2.4]{TW07}.
\begin{lemma}\label{lem:norm_bdd}
    There exists $C=C(\|\phi\|_{L^\infty(I)})>0$ such that
    $$
        \|u(t)\|_{L^\infty(I)}\le C
    $$
    for all $t\in[0,T]$.
\end{lemma}

\begin{proof}
    Let $p>2$.
    Define the operator
    $$
        \mathcal A\varphi
        =
        -d\varphi_{xx}+\varphi,
        \qquad
        D(\mathcal A)=W^{2,p}_{\mathrm{per}}(I).
    $$
    We have $\sigma(\mathcal A)\subset[1,\infty)$.
    Moreover, by \cite{TW07}, $-\mathcal A$ generates an analytic semigroup on $L^p(I)$ satisfying
    \begin{align*}
        \|\mathcal A^\theta e^{-t\mathcal A}v\|_{L^p(I)}
        &\le
        Ct^{-\theta}e^{-\nu t}\|v\|_{L^p(I)},\\
        \|\mathcal A^\theta e^{-t\mathcal A}\partial_xv\|_{L^p(I)}
        &\le
        C(\varepsilon)t^{-\theta-1/2-\varepsilon}e^{-\nu t}\|v\|_{L^p(I)}.
    \end{align*}

    
    We rewrite equation~\eqref{eq:main} as
    \begin{align*}
        u_t
        =
        du_{xx}-u-k\partial_x\bigl[u(K'*u)\bigr]+\widetilde f(u),
    \end{align*}
    where $\widetilde f(u):=f(u)+u$.
    Then,
    \beaa
        u(t)
        &=&
        e^{-t\mathcal A}\phi
        -
        k\int_0^t e^{-(t-s)\mathcal A}
        \partial_x\bigl[u(K'*u)\bigr](s)\,ds+
        \int_0^t e^{-(t-s)\mathcal A}\widetilde f(u(s))\,ds\\
        &=:&u_1(t)+u_2(t)+u_3(t).
    \eeaa

    For $u_1$, we have
    $$
        \|u_1(t)\|_{L^\infty(I)}
        \le
        \|\phi\|_{L^\infty(I)}.
    $$
    Moreover,
    \begin{align*}
        \widetilde f(u)
        =(c+1)u-cu^2
        \le
        \widetilde C
        :=
        \frac{(c+1)^2}{4c}.
    \end{align*}
    Since $e^{-t\mathcal A}$ preserves positivity,
    \begin{align*}
        u_3(t)
        \le
        \widetilde C\int_0^t e^{-(t-s)}\,ds \le \widetilde C.
    \end{align*}


    It remains to estimate $u_2$.
    Since $p>2$, we can choose
    $ \theta\in\left(1/2p,1/2\right)$ and then choose $\varepsilon\in(0,1/2-\theta)$.
    By the embedding $D(\mathcal A^\theta)\hookrightarrow L^\infty(I)$ and the semigroup estimates above,
    \begin{align*}
        \|u_2(t)\|_{L^\infty(I)}
        &\le
        C\int_0^t
        (t-s)^{-\theta-1/2-\varepsilon}e^{-\nu(t-s)}
        \|u(K'*u)(s)\|_{L^p(I)}\,ds.
    \end{align*}
    By Lemma~\ref{lem:Lp_bound} and \eqref{eq:L1_bound},
    \begin{align*}
        \|u(K'*u)\|_{L^p(I)}
        \le
        \|K'_L\|_{L^\infty(I)}
        \|u\|_{L^p(I)}\|u\|_{L^1(I)}
        \le
        C(\|\phi\|_{L^\infty(I)}).
    \end{align*}
    Therefore,
    \begin{align*}
        \|u_2(t)\|_{L^\infty(I)}
        \le
        C(\|\phi\|_{L^\infty(I)})
        \int_0^\infty
        r^{-\theta-1/2-\varepsilon}e^{-\nu r}\,dr
        \le
        C(\|\phi\|_{L^\infty(I)}).
    \end{align*}
    The integral is finite because $\theta+\varepsilon<1/2$.


    By Lemma~\ref{lem:positive}, we have $u\ge0$.
    Since $u=u_1+u_2+u_3$ and $u_3\le\widetilde C$, we obtain pointwise
    $$
        0\le u(t,x)
        \le
        |u_1(t,x)|+|u_2(t,x)|+\widetilde C.
    $$
    Therefore,
    \begin{align*}
        \|u(t)\|_{L^\infty(I)}
        \le
        \|u_1(t)\|_{L^\infty(I)}
        +
        \|u_2(t)\|_{L^\infty(I)}
        +
        \widetilde C.
    \end{align*}
    This proves the conclusion.
\end{proof}

Finally, we give the proof of Proposition~\ref{prop:global_ex}.
 
\begin{proof}[Proof of Proposition~\ref{prop:global_ex}]
    Let $[0,T_{\max})$ be the maximal existence interval of the solution.
    By Lemma~\ref{lem:norm_bdd},
    $$
        \sup_{0\le t<T_{\max}}\|u(t)\|_{L^\infty(I)}
        \le C,
    $$
    where $C$ is independent of $T_{\max}$.
    Set $m:=C$.
    By the proof of Proposition~\ref{prop:local_ex}, there exists $\tau_m>0$ such that every initial value whose $L^\infty$ norm is bounded by $m$ admits a solution on a time interval of length at least $\tau_m$.
    Suppose that $T_{\max}<\infty$.
    Choose $t_0<T_{\max}$ sufficiently close to $T_{\max}$ so that $t_0+\tau_m>T_{\max}$.
    Taking $u(t_0)$ as a new initial value, the local existence result extends the solution to $[t_0,t_0+\tau_m]$.
    By uniqueness, this solution agrees with the original solution on their common interval and therefore extends it beyond $T_{\max}$.
    This contradicts the definition of $T_{\max}$.
    Hence $T_{\max}=\infty$, and equation~\eqref{eq:main} admits a global solution.
\end{proof}

\bibliographystyle{plain}
\bibliography{reference.bib}
\end{document}